\documentclass[preprint,12pt]{elsarticle}

\usepackage[toc,page,title,titletoc,header]{appendix}
\usepackage{stmaryrd}
\SetSymbolFont{stmry}{bold}{U}{stmry}{m}{n}
\usepackage{amsmath}
\usepackage{amssymb}
\usepackage{amsthm}
\usepackage{tabularx}
\usepackage{subfigure}
\usepackage{epsfig}
\usepackage{indentfirst}
\usepackage{cases}
\biboptions{sort&compress}
\usepackage{graphicx}
\usepackage{epsfig}
\usepackage{color}
\usepackage{float}
\usepackage{multirow}
\usepackage{mathtools}

\usepackage{tikz}
\usetikzlibrary{shapes,arrows,positioning}

\usepackage[ruled,vlined]{algorithm2e}
\usepackage{algpseudocode}

\usepackage{caption}

\newcommand{\be}{\begin{equation}}
\newcommand{\ee}{\end{equation}}
\newcommand{\ba}{\begin{array}}
\newcommand{\ea}{\end{array}}
\newcommand{\bea}{\begin{eqnarray}}
\newcommand{\eea}{\end{eqnarray}}
\newcommand{\beas}{\begin{eqnarray*}}
\newcommand{\eeas}{\end{eqnarray*}}

\newtheorem{theorem}{Theorem}[section]

\newtheorem{lemma}[theorem]{Lemma}

\newtheorem{example}[theorem]{Example}

\newcommand{\paren}[1]{\left(#1\right)}
\newcommand{\jump}[1]{\left[#1\right]}

\newcommand{\at}[2]{\left. #1 \right|_{#2}}
\newcommand{\grad}[1]{\nabla #1}
\newcommand{\grads}[1]{\nabla_{\mathcal{S}} #1}
\newcommand{\mb}[1]{\bm{#1}}
\newcommand{\mc}[1]{\mathcal{#1}}

\newcommand{\wt}[1]{\widetilde{#1}}
\newcommand{\bm}[1]{\boldsymbol{#1}}

\newcommand{\vph}{\varphi}

\numberwithin{equation}{section}

\journal{Journal of Computational Physics}

\begin{document}

\begin{frontmatter}



\title{A kernel-free boundary integral method for elliptic interface problems on surfaces}


\author[1]{Pengsong Yin} 

\affiliation[1]{organization={School of Mathematical Sciences, MOE-LSC and Institute of Natural Sciences, Shanghai Jiao Tong University},
            city={Shanghai},
            postcode={200240}, 
            country={China}
}
\ead{yps2701155@sjtu.edu.cn}

\author[1]{Wenjun Ying\corref{4}}
\ead{wying@sjtu.edu.cn}
\cortext[4]{Corresponding author.}

\author[1]{Yulin Zhang}
\ead{yulin.zhang@sjtu.edu.cn}

\author[2]{Han Zhou\corref{4}}
\ead{hzhou24@sas.upenn.edu}

\affiliation[2]{organization={Department of Mathematics, University of Pennsylvania},
            city={Philadelphia},
            postcode={19104}, 
            state={Pennsylvania},
            country={USA}
}

\begin{abstract}
This work presents a generalized boundary integral method for elliptic equations on surfaces, encompassing both boundary value and interface problems. The method is kernel-free, implying that the explicit analytical expression of the kernel function is not required when solving the boundary integral equations. The numerical integration of single- and double-layer potentials or volume integrals at the boundary is replaced by interpolation of the solution to an equivalent interface problem, which is then solved using a fast multigrid solver on Cartesian grids. This paper provides detailed implementation of the second-order version of the kernel-free boundary integral method for elliptic PDEs defined on an embedding surface in $\mathbb{R}^3$ and presents numerical experiments to demonstrate the efficiency and accuracy of the method for both boundary value and interface problems.
\end{abstract}



\begin{keyword}

Kernel-free boundary integral method \sep surface PDE \sep interface problem \sep fast multigrid solver\sep Cartesian grids.



\MSC[2020] 65N06 \sep 35R01 \sep 35J25

\end{keyword}

\end{frontmatter}



\section{Introduction}
Partial differential equations (PDEs) on curved surfaces are crucial for modeling diffusive or transport phenomena on non-planar domains, such as the terrestrial ocean currents~\cite{isern2017remote} or the diffusion of surfactants on fluid interfaces~\cite{langevin2014rheology,muradoglu2014simulations}. Such problems span multiple disciplines including geophysics~\cite{cushman2011introduction,lowrie2020fundamentals}, fluid mechanics~\cite{olshanskii2018finite,reuther2018solving}, and biophysics~\cite{gera2018modeling,elliott2010modeling}. Particularly, an internal sharp interface, a co-dimension two object, often emerges within these surfaces due to the inhomogeneity of the media, exemplified by two-phase flow or phase separation processes. The presence of the sharp interface divides the surface into two or more components, over which the solution is anticipated to be piecewise smooth. Consequently, this paper proposes an efficient numerical method for addressing elliptic PDE interface problems on surfaces.

In recent times, considerable research has been conducted on the development and analysis of numerical methods for solving surface PDEs. Due to the geometric intricacies introduced by the underlying surface, numerical methods are fundamentally categorized as either intrinsic or extrinsic. Intrinsic approaches involve formulating the PDE by globally or locally parameterizing the surface, followed by the discretization of the surface PDE within the parameter space. These approaches encompass surface triangular mesh-based methods, such as surface finite element methods \cite{dziuk2007surface,dziuk2013finite,demlow2012adaptive}, finite volume methods \cite{ju2009finite,calhoun2010finite}, and strategies that involve locally parameterizing the surface by introducing overlapping patches \cite{liang2013solving,fortunato2024high,beale2020solving}. These techniques do not rely on the embedded space, with the computational mesh confined to the surface itself, thereby offering enhanced accuracy and efficiency. Conversely, extrinsic methods involve embedding the surface into an expanded Euclidean space and extending the solution to the PDE beyond the surface. This category includes the closest point method \cite{macdonald2010implicit,macdonald2008level}, the narrow band finite element method \cite{deckelnick2010h,olshanskii2016narrow}, and the trace finite element method \cite{reusken2015analysis,olshanskii2018trace}. The extended solution typically satisfies a simplified PDE, facilitating easier numerical discretization. Moreover, straightforward computational meshes can be generated in the Euclidean space without necessitating alignment with the complex surface.

Concerning solving PDEs with interfaces, the inherent lack of smoothness at sharp interfaces leads to loss of convergence and accuracy for standard discretizations. The jump conditions can be addressed through the employment of an interface-fitted computational mesh. Nonetheless, generating such high-quality fitted triangular meshes poses a non-trivial and time-consuming challenge. Consequently, extensive research over recent decades has concentrated on the development of unfitted mesh-based methods for solving interface problems. This encompasses finite difference-type methods, such as the immersed interface method \cite{li2006immersed,li2001immersed}, the matched interface and boundary method \cite{zhou2006high,yu2007matched}, and the kernel-free boundary integral method (KFBI) \cite{ying2007kernel,ying2014kernel,zhou2024correction}, as well as finite element-type approaches, including the immersed finite element method \cite{zhang2004immersed,lin2015partially}. Unfitted mesh-based methods permit the interface to traverse mesh nodes and devise local strategies to accurately handle the solution jump at the interface while preserving global convergence and efficiency.

There is a paucity of studies concerning numerical methods for surface interface problems, primarily due to the intricate challenges inherent in handling both the underlying surface and the on-surface interface. The local tangential lifting method mitigates this issue by discretizing the surface PDE via a lifting process onto the tangent plane and utilizing an immersed boundary approach with a regularized delta function to solve the moving interface problem on surfaces~\cite{xiao2021local}. Furthermore, there has been a concerted effort to adapt the immersed finite element method for surface interface problems~\cite{guo2023immersed, chen2024mixed, chen2025level}. 

The boundary integral method has demonstrated success in solving interface problems, attributed to its high efficiency through dimension reduction and superior accuracy.
Nonetheless, in the context of a curved two-dimensional manifold, which results in variable coefficients, the absence of a closed-form Green’s function renders the conventional boundary integral method ineffective. Among the recent literature, notable contributions include Goodwill and O'Neil's parametrix method~\cite{goodwill2025parametrix} and Kropinski et al.'s FMM-accelerated Laplace-Beltrami solver~\cite{kropinski2014fast, kropinski2016integral}.

The key advantage of the KFBI method lies in its generalization of traditional boundary integral equation (BIE) frameworks while avoiding their inherent limitations. For example, in interface problems with discontinuous coefficients, the KFBI approach naturally enforces jump conditions via implicit boundary integrals, bypassing the cumbersome layer-potential evaluations required by classical methods. Furthermore, its compatibility with standard PDE discretizations (e.g., finite differences or finite volumes) enables seamless integration into existing computational workflows. 

This work extends the KFBI methodology to elliptic PDEs on curved surfaces, emphasizing its capability to handle both boundary value and interface problems with unmatched simplicity and robustness. Numerical experiments validate the method’s performance against analytical solutions and benchmark problems, underscoring its potential as a versatile tool for surface PDE simulations.

The structure of this paper is organized as follows: Section~\ref{sec:elliptic eq on surf} collects some preliminary results on the geometry of embedding surfaces and formulates the surface elliptic PDE problems, including the boundary value problems and the interface problems. Section~\ref{sec:bif} introduces the boundary integral formulations corresponding to various surface PDEs. Based on the boundary integral formulations, kernel-free boundary integral methods for solving surface boundary value problems and interface problems are proposed in section~\ref{sec:KFBI}. The algorithms are summarized in section~\ref{sec:Algorithm sum}. Extensive numerical experiments are reported in section~\ref{sec:numer}, demonstrating the accuracy and efficiency of the proposed methods. Finally, some conclusions are drawn in section~\ref{sec:conclusion}.

\section{Elliptic PDEs on surfaces}\label{sec:elliptic eq on surf}
\subsection{Geometry of embedding surfaces}

\begin{figure}[htbp]
\centering
\includegraphics[scale=0.4]{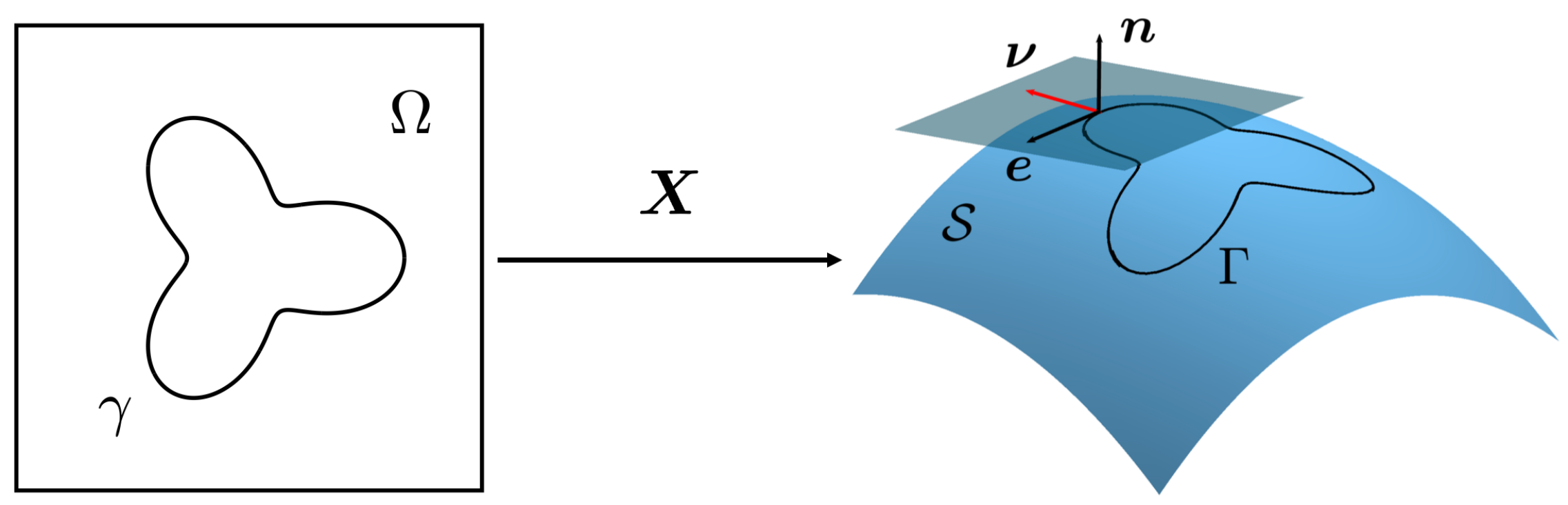}
\caption{A schematic illustration of a parametric surface and an irregular domain defined on it.}
\label{fig:inter problem on surf}
\end{figure}

As illustrated in Fig.~\ref{fig:inter problem on surf}, consider a smooth 2-dimensional Riemannian manifold $\mathcal{S} \subset \mathbb{R}^3$, which can be closed or open. Suppose $\mc S$ has a regular parameterization $\bm X(\bm \xi)$, where $\bm \xi = (\xi_1, \xi_2)\in\mathbb{R}^2$ and $\bm X\colon\mathbb{R}^2\to\mathbb{R}^3$ is the embedding map. We denote by $\Omega\subset\mathbb{R}^2$ the preimage of $\mathcal{S}$, and assume, without losing generality, that $\Omega$ is a planar rectangular domain. Let $\Gamma = \bm X(\gamma)$ be a closed smooth curve on $\mc S$, where $\gamma \subset \Omega$ is a simple closed curve defined by $\bm r(s)=(r_1(s), r_2(s))^T$ with $0<s<|\gamma|$ being the arc-length parameter. Suppose $\gamma$ separates $\Omega$ into the interior region $\Omega^+$ and the exterior region $\Omega^-$ whose images under $\bm X$ are $\mc S^+$ and $\mc S^-$, satisfying $\Gamma = \mc {S^+} \cap \mc S^-$.

Denote $\partial_i \coloneqq \partial_{\xi_i}$ for $i=1,2$. The Riemannian metric tensor on $\mc S$ can be defined as
\begin{equation}
    g_{ij}\coloneqq\partial_i \bm X \cdot\partial_j \bm X,\qquad g\coloneqq\det(g_{ij}),\qquad (g^{ij})\coloneqq(g_{ij})^{-1},\qquad i,j=1,2.
\end{equation} Let $\bm e$ denote the unit tangent vector of $\Gamma$, $\bm n$ the normal vector of $\mc S$, and let $\bm \nu$ be the outer conormal vector, which is considered tangent to $\mc S$ and normal to $\Gamma$, i.e.
\begin{equation}
    \bm e = \partial_s (\bm X\circ\bm r) = \sum_{i=1}^2 \partial_s r_i \partial_i\bm X,\qquad \bm n = \dfrac{\partial_1\bm X\times \partial_2\bm X}{|\partial_1\bm X\times \partial_2\bm X|},\qquad\bm{\nu} = \bm{e}\times\bm{n}.
\end{equation}

The differential operators with respect to $\mc S$ are defined as follows: 
\begin{equation}
\begin{aligned}
    \grads f\coloneqq\sum_{1\leq i,j\leq 2}g^{ij}\partial_{i}f\partial_{j}\bm{X},\qquad\grads \cdot\bm{f}\coloneqq \text{Tr}(\grads \bm{f})=\sum_{i=1}^2\frac{1}{\sqrt{g}}\partial_i(\sqrt{g}f_i),
\end{aligned}
\end{equation}
where $\bm{f}=f_1\partial_1\bm{X}+f_2\partial_2\bm{X}$. 
Consequently, the Laplace-Beltrami operator $\Delta_{\mc S}$ (also termed the surface Laplacian) on $\mc S$ is given by \begin{equation}
    \Delta_{\mc S} f :=\grads \cdot\grads f=\sum_{1\leq i,j\leq 2}\frac{1}{\sqrt{g}}\partial_i\left(\sqrt{g}g^{ij}\partial_{j}f\right).
\end{equation} 

Given a piecewise $C^1$ function $f\colon\mc{S}\to\mathbb{R}$, we define the jump of $f$ and its normal derivative across the interface at $\bm p\in \Gamma$. To this end, consider a curve $\bm y(s)$ on $\mc S$ such that $\bm y(0) = \bm p$ and $\bm y'(0) = \bm \nu(\bm p)$, where $\bm\nu$ is the conormal vector. The jump of $f$ and its normal derivative is then defined as follows:
\begin{subequations}
    \begin{align}
        &\jump{f}(\bm p) = f|_{\Gamma_+}(\bm p) - f|_{\Gamma_-}(\bm p)\coloneqq\lim_{\varepsilon\to 0^+} (f \circ \bm y)(-\varepsilon) - (f\circ \bm y)(\varepsilon),\\
        \begin{split}
            &\jump{\bm \nu \cdot \grads f }(\bm p) = \bm \nu \cdot \grads f |_{\Gamma_+}(\bm p) - \bm \nu \cdot \grads f |_{\Gamma_-}(\bm p) \\
            &\qquad\qquad\quad\,\,\,\coloneqq \lim_{\varepsilon\to 0^+} \bm \nu \cdot \paren{(\grads f \circ \bm y)(-\varepsilon) - (\grads f\circ \bm y)(\varepsilon)}.
        \end{split}
    \end{align}
\end{subequations}

\begin{lemma}[Divergence theorem, \cite{barrett2020parametric}]
    Let $\mc S$ be a compact orientable $C^2$-hypersurface with normal
 vector field $\bm{n}$, let $\bm{f}\in [C^1(\mc S)]^3$. Then it holds that \begin{equation}\label{eqn:div thm}
     \int_{\mc S}\grads \cdot\bm{f}+H\bm{f}\cdot\bm{n}\,\mathrm{d}\mu=\int_{\partial S}\bm{f}\cdot\bm{\nu}\,\mathrm{d}s,
 \end{equation} where $\mu$ is the Riemannian measure on $\mc S$, $H$ denotes the mean curvature and $\bm{\nu}$ represents the outer unit conormal to $\partial \mc S$.
\end{lemma}

\begin{lemma}[Green's identity, \cite{barrett2020parametric}]
    Suppose $f\in C^2(\mc S),g\in C^1(\mc S)$, then \begin{equation}\label{eqn:green id}
        \int_{\mc S}g\Delta_{\mc S}f+\grads f\cdot\grads g\,\mathrm{d}\mu=\int_{\partial S}g\grads f\cdot\bm{\nu}\,\mathrm{d}s.
    \end{equation}
\end{lemma}

\subsection{Boundary value problem}
We consider boundary value problems on an open surface $\mc S^+$ with a smooth boundary $\partial \mc S^+$. Let $\kappa\geq 0$ be a constant coefficient. Given the function $f:\mc S^+\to \mathbb R$, the boundary value problem of an elliptic PDE is formulated as: finding a solution $u:\mc S^+\to \mathbb R$ such that
\begin{equation} \label{boundary value problem}
     \Delta_{\mc S}  u  - \kappa u = f, \qquad  \text{ in }\mc S^+,
\end{equation}
subject to either the Dirichlet boundary condition
\begin{equation}
u = g_D,\qquad\text{ on } \partial\mc S^+ , \label{Dirichlet BVP}
\end{equation}
or the Neumann boundary condition
\begin{align} 
    \bm \nu \cdot \grads u &= g_N,\qquad\text{ on } \partial\mc S^+ .\label{Neumann BVP}
\end{align}
Here, $g_D$ and $g_N$ are the Dirichlet and Neumann boundary data, respectively, and $\bm\nu$ is the outer unit conormal vector to $\partial\mc S^+$.  
The Dirichlet problem admits a unique solution for all $\kappa\ge0$, whereas the Neumann problem has a null space consisting of constant functions and is uniquely solvable only when $\kappa>0$. If $\kappa=0$, the Neumann data $g_N$ must satisfy the compatibility condition
\begin{equation}
\int_{\mc S^+} f \, \mathrm{d}\mu = \int_{\partial\mc S^+} g_N \, \mathrm{d}s,
\end{equation}
in which case the solution is determined only up to an additive constant.

\subsection{Interface problem} 
For the elliptic interface problem on the surface $\mc S$. Let $\beta(\bm p)$ and $\kappa(\bm p)$ be piecewise constant diffusion and reaction coefficients, with $\beta\equiv\beta^\pm>0$ and $\kappa\equiv\kappa^\pm\ge0$ in $\mc S^\pm$. For an open surface, we further assume that the internal interface $\Gamma$ lies away from the boundary $\partial\mc S$; that is, $\mathrm{dist}(\Gamma,\partial\mc S)>0$. Suppose $f$ is a piecewise smooth function on $\mc S$. Then the elliptic interface problem for the solution $u$ is stated as follows:
\begin{equation} \label{Interface problem}
    \grads\cdot\left(\beta\grads u\right)-\kappa u = f, \qquad \text{in } \mc S\setminus\Gamma,
\end{equation} 
subject to interface conditions
\begin{subequations}\label{eqn:inter cond}
    \begin{align}
        \left[u\right]&=g_1,\qquad \text{ on }\Gamma,\label{interface condition1}  \\
        \left[\beta\bm \nu \cdot \grads u\right]&=g_2,\qquad \text{ on }\Gamma,\label{interface condition2} 
    \end{align}
\end{subequations}
where $g_1$ and $g_2$ are two given functions defined on $\Gamma$.
If $\mc S$ is open, one may also impose an additional condition on $\partial\mc S$, such as a Dirichlet boundary condition
\begin{equation}
u = g_0 \qquad \text{on }\partial\mc S,
\end{equation} 
where $g_0$ is the prescribed boundary data.  
If $\mc S$ is closed (and hence has no boundary), then when $\kappa \equiv 0$ the problem has a nontrivial nullspace.  
Therefore, in such cases we assume that either $\kappa^+\neq 0$ or $\kappa^-\neq 0$, so that the problem admits a unique solution.

\section{Boundary integral formulations}\label{sec:bif}
The boundary value and interface problems described in the previous section are reformulated as boundary integral equations in this section.
\subsection{Integral operators}\label{sec:3.1}
Suppose $\bm p$ and $\bm q$ be points in $\mc S$.  Let $\delta_{\bm p}(\bm q)$ be the surface Dirac delta function at $\bm p$, defined as a distribution on $\mc S$ satisfying
\begin{equation}
\int_{\mc S} \phi(\bm q)\,\delta_{\bm p}(\bm q)\,\mathrm{d}\mu
= \phi(\bm p),
\qquad
\forall\,\phi\in C_c^\infty(\mc S).
\end{equation}
We define the Green’s function $G(\bm p,\bm q)$ on $\mc S$ with respect to the elliptic operator $\Delta_{\mc S}-\kappa$ by
\begin{equation}
(\Delta_{\mc S,\bm q}-\kappa)\,G(\bm p,\bm q)
=\delta_{\bm p}(\bm q),
\qquad \bm q\in\mc S,
\end{equation}
where $\Delta_{\mc S,\bm q}$ is the Laplace-Beltrami operator acting on the $\bm q$-variable.  For an open surface $\mc S$, the following homogeneous boundary condition is imposed to guarantee its existence and uniqueness:
\begin{equation}
G(\bm p,\bm q)=0,
\qquad \bm q\in\partial\mc S.
\end{equation}

The Laplace-Beltrami operator can be viewed as a variable‑coefficient elliptic differential operator. In certain special geometries, such as the sphere, one can write down the Green’s function explicitly \cite{kropinski2016integral}. However, for a generic smooth surface, no closed‑form expression is available. Rather than resorting to parametrix or Levi‐function constructions \cite{goodwill2025parametrix, pan2024levi}, which trade off the general utility of the Green’s function for analytic formulas, we rely solely on the layer‐potential framework built from the Green’s function’s abstract properties, without requiring its explicit form.

For a point $\bm p\in\mc S$, we define the following potential functions,
\begin{subequations}\label{eqn:pot}
    \begin{align}
        \text{Yukawa potential}\colon (Vf)(\bm p) &\coloneqq\int_{\mc S}G(\bm{p},\bm{q})f(\bm{q})\,\mathrm{d}\mu_{\bm{q}},\label{eqn:yukawa pot}\\
        \text{Single-layer potential}\colon (S\psi)(\bm p) &\coloneqq \int_{\Gamma}\psi(\bm{q})G(\bm{p},\bm{q})\,\mathrm{d}s_{\bm{q}},\label{eqn:singl pot}\\
        \text{Double-layer potential}\colon (D\varphi)(\bm p) &\coloneqq \int_{\Gamma}\varphi(\bm{q})\bm \nu(\bm{q})\cdot\nabla_{\mc S,\bm{q}}G(\bm{p},\bm{q}) \,\mathrm{d}s_{\bm{q}}.\label{eqn:doubl pot}
    \end{align}
\end{subequations}

For $\bm p\in \mc S\setminus\Gamma$, the integrands are smooth, so the corresponding potentials are finite and smooth on $\mc S\setminus\Gamma$.  Moreover, acting the differential operator $\Delta_{\mc S}-\kappa$ to the potential functions \eqref{eqn:pot} yields \begin{equation}
    (\Delta_{\mc S}-\kappa)\,Vf = f,\qquad (\Delta_{\mc S}-\kappa)\,S\psi = 0,\qquad (\Delta_{\mc S}-\kappa)\,D\varphi = 0, \qquad\text{in}\,\,\mc{S}\backslash\Gamma.
\end{equation}

For $\bm p\in\Gamma$, both $V f$ and $S\psi$ are continuous at $\bm p$. Hence, traces of $Vf$ and $S\psi$ on $\Gamma$ are denoted using the same notation for simplicity. Since the integral $D\vph$ on $\Gamma$ is divergent, we interpret it as a principal value integral and denote it by $K\varphi$, i.e.
\begin{equation}
    (K\varphi)(\bm p) \coloneqq \text{p.v}.\int_{\Gamma}\varphi(\bm{q})\bm \nu (\bm{q})\cdot\nabla_{\mc S,\bm{q}}G(\bm{p},\bm{q}) \,\mathrm{d}s_{\bm{q}} = \frac{1}{2}\paren{(D\varphi)|_{\Gamma_+} + (D\varphi)|_{\Gamma_-}}(\bm p).
\end{equation}
Correspondingly, the double-layer potential is discontinuous on $\Gamma$ and satisfies the following jump relation,
\begin{align}
    (D\vph)|_{\Gamma^+}(\bm p) - (D\vph)|_{\Gamma^-}(\bm p) = \vph(\bm p).
\end{align}

The normal derivative of $V f$ remains continuous across $\Gamma$, whereas the normal derivative of the single-layer potential satisfies the following jump relation
\begin{align}
    \bm \nu \cdot \grad_\mathcal{S} (S\psi)|_{\Gamma^+}(\bm p) - \bm \nu \cdot \grad_\mathcal{S} (S\psi)|_{\Gamma^-}(\bm p) = -\psi(\bm p).
\end{align} The adjoint double-layer integral is defined by
\begin{equation}
\begin{aligned}
(K^\prime\varphi)(\bm p) &\coloneqq \text{p.v.}\int_{\Gamma}\varphi(\bm{q})\bm \nu(\bm{p})\cdot\nabla_{S,\bm{p}}G(\bm{p},\bm{q})\,\mathrm{d}s_{\bm{q}} \\
&= \frac{1}{2}\paren{\bm \nu \cdot \grad_\mathcal{S} (S\psi)|_{\Gamma^+}+\bm \nu \cdot \grad_\mathcal{S} (S\psi)|_{\Gamma^-}}(\bm p).
\end{aligned}
\end{equation}and its normal derivative is also continuous. The hyper-singular integral operator is defined as follows
\begin{align}
(H\varphi)(\bm p)\coloneqq\text{p.f.}\int_{\Gamma}\bm \nu^T(\bm{p})\left(\nabla_{S,\bm{p}}\nabla_{\mc S,\bm{q}}G(\bm{p},\bm{q})\right)\bm \nu(\bm{q})\varphi(\bm{q})\,\mathrm{d}s_{\bm{q}},
\end{align}
where the notation $\text{p.f.}$ means the Hadamard finite part integral. 

For a detailed discussion of the above integral operators and their associated jump relations, we refer the reader to \cite{mitrea2000potential,mitrea1999boundary}.

\subsection{Boundary value problem} \label{subsection: Boundary value problem}

Analogous to the planar case in \cite{ying2007kernel}, the solution $u(\bm p)$ to the Dirichlet BVP \eqref{boundary value problem}--\eqref{Dirichlet BVP} can be expressed as a sum of a Yukawa potential and a double-layer potential: 
\begin{equation} \label{solution for Dirichlet}
    u(\bm p)= Vf(\bm p)+ D\varphi(\bm p),\qquad\text{in } \mc S,
\end{equation} 
where the density $\varphi(\bm{p})$ satisfies a Fredholm boundary integral equation (BIE) of the second kind as follows: 
\begin{equation} \label{BIE for Dirichlet}
\frac{1}{2}\varphi(\bm p)+ K\varphi(\bm p)+ Vf(\bm p)=g_D(\bm p),\qquad\text{on } \Gamma .
\end{equation} 
The derivation follows an approach similar to that in \cite{ying2014kernel} by applying the Green's identities \eqref{eqn:div thm}--\eqref{eqn:green id}. Details are omitted here for brevity.


Similarly, the solution to the Neumann BVP \eqref{boundary value problem}--\eqref{Neumann BVP} can be expressed in terms of a Yukawa potential and a single-layer potential: \begin{equation} \label{solution for Neumann}
    u(\bm p)= Vf(\bm p)- S\psi(\bm p), \qquad\text{in } \mc S.
\end{equation} And the density function $\psi(\bm p)$ is required to satisfy the following BIE: \begin{equation} \label{BIE for Neumann}
    \frac{1}{2}\psi(\bm p)- K^\prime\psi(\bm p)+\bm \nu \cdot\grads (Vf)(\bm p)=g_N(\bm p),\qquad\text{on } \Gamma .
\end{equation}

Since the Green's function can be viewed as a perturbation of the Green's function in Euclidean space by a more regular term, the operators $ K $ and $ K' $ remain compact, with spectra clustered around zero. Consequently, the operators $ \frac{1}{2}\mathcal{I} + K $ and $ \frac{1}{2}\mathcal{I} - K' $ have spectra clustered around $ 1/2 $, which is away from zero. This spectral property enables the use of Krylov subspace iterative methods, such as GMRES, to efficiently solve the resulting discrete algebraic system, provided the discrete operators exhibit similar spectral behavior.

\subsection{Interface problem with $\kappa^+/\beta^+ = \kappa^-/\beta^-$}
In this case, one can derive a simple boundary integral equation of the second kind with one unknown density function. To this end, we define the following Green's function:
\begin{equation}
        \paren{ \Delta_{\mc S, \bm{q}} - \frac{\kappa^\pm}{\beta^\pm}}G(\mb p, \bm q)= \delta_{\bm p}(\bm q),  \qquad  \bm q\in\mc S.
\end{equation} 

Let $\psi$ be an unknown single-layer density. In a manner similar to \cite{ying2014kernel}, one can express the solution as a linear combination of the single-layer, double-layer and Yukawa potentials, i.e.
\begin{equation} \label{soltion for interface problem}
    u(\bm p) = D g_1 (\bm p) - S\psi(\bm p) + V  \hat{f} (\bm p), \qquad \hat{f}(\bm p) \coloneqq f(\bm p) /\beta(\bm p), \qquad \bm p\in \mathcal{S}.
\end{equation}

One can easily verify that the above representation satisfies the PDE~\eqref{Interface problem}, the interface condition~\eqref{interface condition1}, and the jump condition $\jump{\bm \nu \cdot \grads u} = \psi$. As long as $u$ also satisfies the interface condition~\eqref{interface condition2}, it will serve as a solution to the interface problem \eqref{Interface problem}--\eqref{eqn:inter cond}.

By taking the derivative of $u$ in the conormal direction, we obtain \begin{equation}
    H g_1  - K^\prime \psi + \bm \nu \cdot \grad_{\mc S}(V \hat{f} )  = 
        \bm \nu \cdot\grad_{\mc S} u|_{\Gamma^+} - \dfrac{1}{2}\psi = \bm \nu \cdot\grad_{\mc S} u|_{\Gamma^-} + \dfrac{1}{2}\psi, \qquad  \bm p\in \Gamma.
\end{equation}
Applying the interface condition \eqref{interface condition2} results in a boundary integral equation of the second kind for the density $\psi$: \begin{equation}\label{eqn:interface-bie1}
    \psi - 2 A_{\beta} K^\prime \psi = \dfrac{2}{\beta^+ + \beta^-} g_2 - 2 A_\beta (H g_1 + \bm \nu \cdot \grad_{\mc S}(V \hat{f} ) ):= \hat{g}_2  ,\qquad\text{ on }\Gamma,
\end{equation}
where $A_{\beta} = (\beta^+ - \beta^-)/(\beta^+ + \beta^-) \in (-1,1)$ is the Atwood ratio.
For $\beta^\pm \neq 0$, the integral equation above with homogeneous right-hand side admits only the zero solution. Hence, by the Fredholm alternative, the boundary integral equation~\eqref{eqn:interface-bie1} is uniquely solvable and has a bounded inverse.

\subsection{Interface problem with $\kappa^+/\beta^+ \neq \kappa^-/\beta^-$}
For the generic case, the Green's functions associated with the differential operators in $\mc S^+$ and $\mc S^-$ must be defined separately. Let $G^{\pm}(\bm{p}, \bm{q})$ be the Green's function associated with the interface problem with variable coefficient in $\mathcal{S}^+$, which satisfies
\begin{equation}
        \paren{ \Delta_{\mc S, \bm{q}} - \frac{\kappa^\pm}{\beta^\pm}}G^\pm(\mb p, \bm q)= \delta_{\bm p}(\bm q),  \qquad  \bm q\in\mc S.
\end{equation}

Similar to the previous work \cite{ying2014kernel}, the single layer, double layer, adjoint double layer, and hyper-singular boundary integrals need to be defined for both the interior and exterior regions. We denote the associated integral operators as $V^\pm,K^\pm, S^\pm, (K^{\prime})^\pm, H^\pm$. 

Let $\varphi=u^+,\psi=\bm{\nu}\cdot\grads u^-$ be two unknown density functions. Applying the interface conditions \eqref{eqn:inter cond}, we can derive the following system of equations: \begin{subequations}
    \begin{align}
        \label{sys subq1} &\frac{1}{2}\varphi=K^+\varphi-S^+\psi-S^+g_N+V^+f^+,\\  
        \label{sys subq2} &\frac{1}{2}(\varphi-g_D)=-K^-\varphi+S^-\psi+K^-g_D+V^-f^-,\\
        \label{sys subq3} &\frac{1}{2}(\psi+g_N)=H^+\varphi-(K^{\prime})^+\psi-(K^{\prime})^+ g_N+\beta^+ \bm \nu \cdot \grads V^+f^+,\\
        \label{sys subq4} &\frac{1}{2}\psi=-H^-\psi+(K^{\prime})^-\psi+H^-g_D+\beta^- \bm \nu \cdot\grads V^-f^-.
    \end{align}
\end{subequations}

Adding \eqref{sys subq1} to \eqref{sys subq2}, \eqref{sys subq3} to \eqref{sys subq4}, and making rearrangement of the terms, the interface problem is reformulated as a system of two boundary integral equations: for $\bm q \in \mc S$, \begin{subequations}
\begin{align}
& \varphi-\left(D^+-D^-\right) \varphi+\left(S^+-S^-\right) \psi=r(\mathbf{q}), \\
& -\left(H^+-H^-\right) \varphi+\psi+\left( (K^{\prime}) ^ + -(K^{\prime}) ^ -\right) \psi=s(\mathbf{q}),
\end{align}
\end{subequations} where $r,s$ are given by \begin{subequations}
    \begin{align}
& r(\mathbf{q})\coloneqq\frac{1}{2} g_D+V^+ f^+ +V^- f^--S^+ g_N+D^- g_D, \\
& s(\mathbf{q})\coloneqq-\frac{1}{2} g_N+\beta ^ +\bm \nu \cdot \grads V^+ f^+ +\beta^-\bm \nu \cdot \grads V^- f^- -(K^{\prime})^+ g_N+H^- g_D.
\end{align}
\end{subequations}

After solving the boundary integral equation system to determine $\varphi$ and $\psi$, the solution $u_i$ and $u_e$ are then obtained by evaluating the following integrals:
\begin{subequations}
\begin{align}
    &u_i(\bm{p}) = D^+ \varphi -S^+\psi-S^+g_N+V^+f^+, \\
    &u_e(\bm{p}) = -D^- \varphi +S^-\psi+S^-g_N+V^-f^-.
\end{align}
\end{subequations}

\section{Kernel-free boundary integral method}\label{sec:KFBI}
In this section, we introduce the kernel-free boundary integral method for evaluating boundary and volume integrals by solving equivalent interface problems on a Cartesian grid with a finite difference method.

\subsection{Equivalent interface problems} \label{Equivalent interface problems}

Rather than interpreting the Yukawa potential, the single-layer potential, and the double-layer potential through their integral representations, we regard them as solutions of simpler interface problems, with their jump relations given in Section~\ref{sec:3.1}.

The interface problems corresponding to the three potentials $ Vf $, $ S\psi $, and $ D\vph $ can be unified in the following form, \begin{subequations}\label{eqn:equiv inter prob}
\begin{align}
    \Delta_{\mc S} u - \kappa u &= F, \qquad \text{in } \mc S \setminus \Gamma, \\
    \jump{u} &= \Phi, \qquad \text{on } \Gamma, \\
    \jump{\bm \nu \cdot \grads u} &= \Psi, \qquad \text{on } \Gamma,
\end{align}
\end{subequations} where the source term $ F $ and jump data $ \Phi $, $ \Psi $ are specified as follows,
\begin{itemize}
    \item $ u = Vf $: $ F = f $, $ \Phi = 0 $, $ \Psi = 0 $;
    \item $ u = S\psi $: $ F = 0 $, $ \Phi = 0 $, $ \Psi = -\psi $;
    \item $ u = D\vph $: $ F = 0 $, $ \Phi = \vph $, $ \Psi = 0 $.
\end{itemize}

These interface problems are significantly simpler than the boundary value problem~\eqref{boundary value problem} and the original interface problem~\eqref{Interface problem}, as they involve constant coefficients and known density functions. This enables the construction of numerical schemes that require only right-hand side corrections and can be efficiently solved using fast elliptic solvers. As a result, one no longer needs an explicit analytical expression for the Green's function to perform numerical quadrature, and (nearly) singular integrals can be completely avoided.




For numerical convenience, we pullback functions on $\mc S$ to the planar domain $ \Omega$. For a surface function $f\colon\mc{S}\to\mathbb{R}$, we denote its pullback under the parametrization $\bm X$ by $\wt{f} = f\circ \bm X$. Then the interface problem \eqref{eqn:equiv inter prob} can be reformulated as an interface problem on a planer domain with variable coefficients
\begin{subequations}\label{eqn:var-ifp}
    \begin{align}
        \sum_{i,j=1}^2 \partial_i\paren{ a_{ij}(\bm x)\partial_j\wt{u}} (\bm x)- a(\bm x) \wt{u}(\bm x)&= \wt{f}(\bm x), \qquad \text{ in }\Omega\setminus\gamma,\\
        \jump{\wt u}(\bm x) &= \wt \Phi(\bm x),\qquad \text{ on }\gamma,\\
        \sum_{i,j=1}^2 \beta_{ij}(\bm x)\jump{ \partial_i \wt u} (\bm x)&= \wt \Psi(\bm x), \qquad \text{ on }\gamma.
    \end{align}
\end{subequations}
where the spatially dependent functions are given by
\begin{equation}
    a=\kappa \sqrt{g}, \qquad \beta_{ij}=g^{ij}\paren{\wt {\bm \nu} \cdot\partial_j\bm X}, \qquad a_{ij} = \sqrt{g} g^{ij}, \qquad f = \sqrt{g} F.
\end{equation}

\subsection{Corrected finite difference scheme}\label{Corrected finite difference scheme}
A finite difference scheme is devised for solving the interface problem~\eqref{eqn:var-ifp}, with right-hand side correction to incorporate the jump of the solution at the interface.
Consider a variable coefficient self-adjoint differential operator $L$, defined as
\begin{equation}
    L = \sum_{i,j=1}^2 \frac{\partial}{\partial x_i} \paren{a_{ij}(\bm x)\frac{\partial}{\partial x_j}} - a(\bm x),\qquad a_{ij} (\bm x) = a_{ji} (\bm x),\qquad \bm x\in\Omega.
\end{equation}
where $a_{ij}, a$ are spatially dependent functions. Assume that $a_{ij}\in C^1(\Omega)$ is pointwise positive definite and $a\ge 0$ such that the differential operator $L$ is negative definite with a suitable boundary condition on $\partial\Omega$. 

For simplicity, we assume that the planar domain $\Omega$ is a square $\Omega=(-1,1)^2$. Let $N$ be a positive integer. We uniformly partition $\Omega$ into a Cartesian grid in each spatial direction with grid spacing $h = \frac{2}{N}$.
The set of grid nodes are denoted as $\overline{\Omega_h}$, defined as
\begin{equation}
    \overline{\Omega_h} =\left\{\left(x_i, y_j\right)=(-1+i h, -1+j h)\mid i, j=0,1, \cdots, N\right\}.
\end{equation}
The sets of interior grid nodes $\Omega_h$ and boundary grid nodes $\partial\Omega_h$ are defined as
\begin{equation}\Omega_h =\left\{\left(x_i, y_j\right)=(-1+i h, -1+j h)\mid i, j=1,2, \cdots, N-1\right\}, \qquad \partial \Omega_h=\overline{\Omega_h} \backslash \Omega_h .
\end{equation}

The diagonal terms $\frac{\partial}{\partial x} a_{11}(x, y) \frac{\partial}{\partial x}+\frac{\partial}{\partial y} a_{22}(x, y) \frac{\partial}{\partial y} - a(x,y)$ at a grid node $\bm x=(x,y)\in\Omega_h$ can be discretized with a standard five-point difference scheme, represented as
\begin{equation}\label{eqn:fds-diag}
\begin{aligned}
    &    h^{-2} 
    \begin{bmatrix}
        0 & a_{22}(x,y+h/2) & 0 \\
        a_{11}(x-h/2,y) &
        \left\{\begin{aligned} &-a_{11}(x-h/2,y)-a_{11}(x+h/2,y) \\ 
        &- a_{22}(x,y+h/2) - a_{22}(x,y-h/2)
        \end{aligned}\right\} 
        & a_{11}(x+h/2,y)\\
        0 & a_{22}(x,y-h/2) & 0
    \end{bmatrix} \\
    &\quad- 
    \begin{bmatrix}
        0 & 0 & 0\\
        0 & a(x,y) & 0 \\
        0 & 0 & 0
    \end{bmatrix}.
\end{aligned}
\end{equation}



The discretization of the mixed term $\frac{\partial}{\partial x} a_{12}(x, y) \frac{\partial}{\partial y}+\frac{\partial}{\partial y} a_{12}(x, y) \frac{\partial}{\partial x}$ results in a seven-point difference scheme, whose stencil depends on the sign of the off-diagonal coefficient $a_{12}$. 
If $a_{12}(x,y)=a_{21}(x,y) < 0$, the seven-point difference scheme reads
\begin{equation}\label{eqn:fds-mixed-1}
\renewcommand{\arraystretch}{1.5}
\frac{h^{-2}}{2}\left[\begin{array}{ccc}
-a_{12}^A-a_{12}^B & a_{12}^A+a_{12}^C & 0 \\
a_{12}^B+a_{12}^D & -a_{12}^C-a_{12}^D-a_{12}^E-a_{12}^F & a_{12}^H+a_{12}^F \\
0 & a_{12}^E+a_{12}^G & -a_{12}^H-a_{12}^G
\end{array}\right],
\end{equation} 
where the upper index indicates the evaluation at  $\bm{x}+\delta h \text { with } \delta$  specified below
\begin{equation}
\renewcommand{\arraystretch}{1.5}
\begin{array}{c|cccccccc} 
& \text { A } & \text { B } & \text { C } & \text { D } & \text { E } & \text { F } & \text { G } & \text { H } \\
\hline \delta & \left(-\frac{1}{2}, 1\right) & \left(-1, \frac{1}{2}\right) & \left(0, \frac{1}{2}\right) & \left(-\frac{1}{2}, 0\right) & \left(0,-\frac{1}{2}\right) & \left(\frac{1}{2}, 0\right) & \left(\frac{1}{2},-1\right) & \left(1,-\frac{1}{2}\right).
\end{array}
\end{equation}
If $a_{12}(x,y)=a_{21}(x,y) > 0$, then the seven-point difference scheme is given by 
\begin{equation}\label{eqn:fds-mixed-2}
\renewcommand{\arraystretch}{1.5}
\frac{h^{-2}}{2}\left[\begin{array}{ccc}
0 & -a_{12}^A-a_{12}^C & a_{12}^A + a_{12}^H \\
-a_{12}^B-a_{12}^D & a_{12}^C+a_{12}^D+a_{12}^E+a_{12}^F & -a_{12}^H-a_{12}^F \\
a_{12}^B+a_{12}^G & -a_{12}^E-a_{12}^G & 0
\end{array}\right]
\end{equation}
where the upper index indicates the evaluation at  $\bm{x}+\delta h \text { with } \delta$  specified below: 
\begin{equation}
\renewcommand{\arraystretch}{1.5}
\begin{array}{c|cccccccc} 
& \text { A } & \text { B } & \text { C } & \text { D } & \text { E } & \text { F } & \text { G } & \text { H } \\
\hline \delta & \left(\frac{1}{2}, 1\right) & \left(-1, -\frac{1}{2}\right) & \left(0, \frac{1}{2}\right) & \left(-\frac{1}{2}, 0\right) & \left(0,-\frac{1}{2}\right) & \left(\frac{1}{2}, 0\right) & \left(-\frac{1}{2},-1\right) & \left(1,\frac{1}{2}\right).
\end{array}
\end{equation}

We denote by the associated coefficient matrix of the seven-point scheme as $L_h$. The stencil selection procedure according to the off-diagonal coefficient is mainly for ensuring the stability of the matrix $L_h$ and the solvability of the resulting linear system. In particular, under mild assumptions, the seven-point finite difference scheme is second-order and the negative of $L_h$ is an M-matrix with guaranteed invertibility and positivity \cite[p.~104]{hackbusch2017elliptic}:
\begin{lemma}[\cite{hackbusch2017elliptic}]
    Assume that $a_{ij}\in C^{2,1}(\overline{\Omega})$ and $|a_{12}(x,y)| < \min(|a_{11}(x,y)|, |a_{22}(x,y)|)$, the finite difference scheme~\eqref{eqn:fds-diag},\eqref{eqn:fds-mixed-1}, and \eqref{eqn:fds-mixed-2} has second-order consistency and for sufficiently small $h$ the associated matrix $-L_h$ is a symmetric and positive-definite M-matrix.
\end{lemma}

We apply the seven-point finite difference scheme to the differential equation $L u = f$ whose solution $u$ is only expected to be piecewise smooth in $\Omega^\pm$.
At grid nodes away from the interface where the solution is smooth, the scheme retain its second order consistency.
At grid nodes in the vicinity of the interface, the finite difference scheme may be taken across the interface, leading to large local truncation error due to the non-smoothness of the solution. Such an issue can be addressed by incorporating the leading local truncation error in the finite difference scheme as correction terms to the right-hand side, modifying the local truncation error and leaving the coefficient matrix intact.

Let $u_{ij}$ be the numerical approximation of $u\left(x_i, y_j\right)$. We write the seven-point finite difference scheme in a generic form,
\begin{equation}
L_h u_{ij}=\sum_{\left(x_r, y_s\right) \in \overline{\Omega_h}} c_{ijrs} u_{rs} = \sum_{\left(x_r, y_s\right) \in \mathcal{S}_{ij}} c_{ijrs} u_{rs}, \qquad\left(x_i, y_j\right) \in \Omega_h,
\end{equation} 
where $c_{ijrs}$ is the coefficient depending on $a_{ij}, a$ and $\mc S_{ij}$ is the set of stencil nodes of the finite difference operator $L_h$ at $(x_i,y_j)$, defined as 
\begin{equation}
\mathcal{S}_{ij}=\left\{\left(x_r, y_s\right) \in \overline{\Omega_h} \mid c_{ijrs} \neq 0\right\}.
\end{equation}

We define the sets of the regular and irregular nodes, denote by $\mc{R}_h,\mc{I}_h$,\begin{equation}
\mathcal{R}_h \coloneqq\left\{\left(x_i, y_j\right) \in \Omega_h \mid \mathcal{S}_{ij} \cap \Omega^{+}=\emptyset \text { or } \mathcal{S}_{ij} \cap \Omega^{-}=\emptyset\right\},\qquad\mathcal{I}_h  \coloneqq\Omega_h \backslash \mathcal{R}_h.
\end{equation}
Let $\Omega_\gamma\supset \mc I_h$ be a narrow band near $\gamma$ with sufficient width to cover all irregular nodes. Suppose the piecewise smooth solution $u$ can be smoothly extended into $\Omega^\mp$ with the restricted values $\at{u}{\Omega^\pm}$. denote by the extended smooth functions as $u^\pm$, defined in $\Omega^\pm\cup\Omega_\gamma$. Let $C = u^+- u^-$ be the difference of the two smooth functions in $\Omega_\gamma$, which is referred to as the correction function.
The local truncation error of the finite difference scheme at a grid node $(x_i, y_j)\in\Omega_h$ is given by
\begin{equation}
E_{ij}  = L_h u\left(x_i, y_j\right)-f\left(x_i, y_j\right) = D(x_i, y_j; C) + \mathcal{O} (h^2).
\end{equation}
where the leading term $D(x_i, y_j; C)$ is given by
\begin{equation}
    D(x_i, y_j; C) = 
    \begin{cases}
0, &\left(x_i, y_j\right) \in \mathcal{R}_h, \\
-\displaystyle\sum_{\left(x_r, y_s\right) \in \mathcal{S}_{ij} \cap \Omega^{-}} c_{ijrs} C\left(x_r, y_s\right), &\left(x_i, y_j\right) \in \mathcal{I}_h \cap \Omega^{+}, \\
 \displaystyle\sum_{\left(x_r, y_s\right) \in \mathcal{S}_{ij} \cap \Omega^{+}} c_{ijrs} C\left(x_r, y_s\right), &\left(x_i, y_j\right) \in \mathcal{I}_h \cap \Omega^{-} .
\end{cases}
\end{equation}
Since the coefficient satisfies $c_{ijrs} = \mc O(h^{-2})$, the local truncation error at irregular nodes behaves like $E_{ij} = \mc O(h^{-2})$. An immediate remedy is to use the value of $D(x_i, y_j; C)$ as correction terms to design a corrected scheme.

Since $D(x_i, y_j; C)$ is a linear combination of the correction function $C$, which is generally not known explicitly but can be approximated by a numerical solution $C_h$.
Denote $\mathcal{Z}_h$ as grid nodes where the correction function is used, i.e.
\begin{equation}
\mathcal{Z}_h=\left(\bigcup_{\left(x_i, y_j\right) \in \mathcal{I}_h \cap \Omega^{+}}\left(\mathcal{S}_{ij} \cap \Omega^{-}\right)\right) \bigcup\left(\bigcup_{\left(x_i, y_j\right) \in \mathcal{I}_h \cap \Omega^{-}}\left(\mathcal{S}_{ij} \cap \Omega^{+}\right)\right).
\end{equation} 
Suppose the numerical solution $C_h$ has pointwise third-order accuracy, i.e., $\left|C-C_h\right|\left(x_i, y_j\right)=\mathcal{O}(h^3)$ for $\left(x_i, y_j\right) \in \mathcal{Z}_h$. Then the corrected scheme is given by
\begin{equation}\label{correction of scheme}
\begin{aligned}
L_h u_{ij} & =f\left(x_i, y_j\right)+D(x_i,y_j;C_h), \qquad\left(x_i, y_j\right) \in \Omega_h .
\end{aligned}
\end{equation}
After correction, the local truncation error is modified to
\begin{equation}
\begin{aligned}
    \widetilde{E}_{ij} &= L_h u\left(x_i, y_j\right) - f\left(x_i, y_j\right) - D(x_i,y_j;C_h) \\
    &= E_{ij}- D(x_i,y_j;C_h)= 
\begin{cases}
\mathcal{O}(h^2),  &\left(x_i, y_j\right) \in \mathcal{R}_h, \\ 
\mathcal{O}(h), & \left(x_i, y_j\right) \in \mathcal{I}_h.
\end{cases}
\end{aligned}
\end{equation}
Since irregular nodes only appears in the vicinity of the interface, which is a co-dimension one object, the first-order local truncation error at irregular nodes is sufficient to result in a uniformly second-order accuracy \cite{beale2007accuracy}.
\begin{theorem}
    Suppose the local truncation error satisfies $|\wt E_{ij}| \leq \bar{C}h$ for $(x_i,y_j)\in\mc I_h$ and $|\wt E_{ij} |\leq \bar{C}h^2$ for $(x_i,y_j)\in\mc R_h$. Then we have the error estimate,
    \begin{equation}
        |u(x_i,y_j) - u_{ij}| \leq \bar{C}h^2, \qquad (x_i,y_j)\in\Omega_h,
    \end{equation}
    where $\bar{C}$ is dependent on $u$ but independent of $h$.
\end{theorem}

The corrected finite‐difference scheme leaves the coefficient matrix intact, which permits the use of fast elliptic solvers for the resulting linear system.  In particular, exploiting the nested structure of the Cartesian grid, we employ a geometric multigrid solver with a full‐multigrid cycle~\cite{trottenberg2001multigrid}, achieving optimal complexity $\mathcal O\bigl(N^2\log\varepsilon\bigr)$, where $\varepsilon$ is the prescribed tolerance and $N^2$ is the total number of grid nodes.

After solving the interface problem, we obtain the grid values of the potential functions. However, in solving boundary integral equations, one would need the boundary integrals on the interface, which are either the trace or normal derivatives of the potential functions on the interface. We apply Lagrange interpolation to extract boundary information from the grid values to complete the integral evaluation on the interface. Since potential functions are not smooth on the interface, corrections are also needed to obtain accurate interpolation. We refer the reader to \cite{ying2007kernel, zhou2024correction} for details of interpolation with the correction function.


\subsection{Local Cauchy problem}\label{Local Cauchy problem}
Here we introduce numerical approach to obtain the correction function $C_h$.
Note that there are infinite many ways to smoothly extend a smooth function. Here we use the approach based on extending the original PDE.
Suppose the piece-wise smooth right-hand side function $f$ can be smoothly extended into $\Omega_\gamma$ from each side to give two smooth function $f^+$ an $f^-$. The extension of $u$ is done by letting $u^\pm$ be solutions to the equation $Lu^\pm = f^\pm$ in the domain $\Omega^\pm\cap \Omega_\gamma$. By elliptic regularity theory, $u^\pm$ and their difference, the correction function, are expected to be smooth.  
Denote by $\bar{f} = f^+ - f^-$.
It is straightforward to show that the correction function satisfies the Cauchy problem of an elliptic PDE as follows
\begin{subequations}
    \begin{align}
        \sum_{i,j=1}^2 \partial_i\paren{ a_{ij}\partial_j C} - a C&= \bar{f}, \qquad \text{ in }\Omega_\gamma,\\
        C &= \bar{\Phi},\qquad \text{ on }\gamma,\\
        \sum_{i,j=1}^2 \beta_{ij}\partial_i C &= \bar{\Psi}, \qquad \text{ on }\gamma.
    \end{align}
\end{subequations}
where $\bar{\Phi}, \bar{\Psi}$ are smooth Cauchy data given on $\gamma$.
Since we only need the value of $C$ at irregular nodes $\mc I_h\subset \Omega_\gamma$, we only solve Cauchy problem for a short distance away from $\gamma$.

Let $\{\bm q_l\}_{l=1}^m\subset \gamma$ be a set of quasi-uniform distributed points on $\gamma$.  Let $B_r(\bm q_l)$ be a ball centered at $\bm q_l$ with radius $r$. We set $r$ sufficiently large such that $\cup_{l=1}^M B_r(\bm q_l) \supset \Omega_\gamma$ and, typically, $r = \mc O(h)$.
The numerical solution is represented using a partition of unity
\begin{equation}
    C_h(\bm x) = \sum_{l = 1}^M \omega_l(\bm x) C_{h, l}(\bm x),\qquad \bm x\in \cup_{l=1}^M B_r(\bm q_l).
\end{equation}
where $\omega_l$ is an arbitrary partition of unity function subordinate to the ball $B_r(\bm q_l)$ and $C_{h, l}$ is the local solution in $B_r(\bm q_l)$.
The simplest partition of unity is by defining it as
\begin{equation}
    \omega_l(\bm x) = 
    \begin{cases}
        1, & l = \operatorname*{arginf}\limits_{k} |\bm q_k -\bm x| ,\\
        0, & \text{else},
    \end{cases}\qquad \bm x\in\Omega_\gamma.
\end{equation}
Such a partition of unity function is piecewise constant with value either $0$ or $1$, so that
\begin{equation}
    C_h(\bm x) = C_{h, l}(\bm x),\qquad \text{if } l = \operatorname*{arginf}_{k} |\bm q_k -\bm x|.
\end{equation}
In such way, we can separately solve the local solutions $C_{h,l}$ in each domain $B_r(\bm q_l)$. The local solution $C_{h, l}$ is represented as a linear combination of polynomial basis
\begin{equation}
    C_{h, l}(\bm x) = \sum_{k = 1}^6 \alpha_{l, k} \phi_k\paren{r^{-1}(\bm x - \bm q_l)},\qquad \bm x\in B_r(\bm q_l),
\end{equation}
where $\alpha_{l,k}$ is the undetermined coefficients and $\phi_k$ is the monomial basis
\begin{equation}
\{\phi_j(x,y)\}_{j=1}^6 = \{1,x,y,x^2,y^2, xy\}.
\end{equation}
The prefactor $r^{-1}$ is designed to spatially rescale the local problem so the argument of $\phi_k$ has size $\mc O(1)$ rather than $\mc O(r)$.

We seek the value of $\alpha_{l,k}$ by collocation.
We locally parameterize the curve segment $\gamma\cap B_r(\bm q_l)$ by $\bm r(s)$ for $s\in [-1,1]$ with $\bm r(0) = \bm q_l$.
Generally, one can select many collocation points from $B_r(\bm q_l)$ and $\gamma\cap B_r(\bm q_l)$ for the PDE and boundary conditions, to result in an overdetermined system and solve in the least squares sense.
Here, to minimize computational cost, we will select the minimum number of collocation points to obtain a squared linear system based on \cite{zhou2024correction}.
We select collocation points $\{\bm p_{l,s}^{(1)}\}_{s=1}^3$ for the Dirichlet boundary condition, $\{\bm p_{l,s}^{(2)}\}_{s=1}^2$ for the Neumann boundary condition, and $\{\bm p_{l}^{(3)}\}$ for the PDE as
\begin{equation}
    \begin{aligned}
        &\bm p_{l,1}^{(1)} = \bm r(0),\quad \bm p_{l,2}^{(1)} = \bm r(-1),\quad \bm p_{l,3}^{(1)} = \bm r(1),\\
        &\bm p_{l,1}^{(2)} = \bm r(-1),\quad \bm p_{l,2}^{(2)} = \bm r(1),\quad \bm p_{l}^{(3)} = \bm q_l.
    \end{aligned}
\end{equation}
Now we can form the discrete collocation problem as \begin{subequations}\label{eqn:dis-cau}
    \begin{align}
        \sum_{i,j=1}^2 \partial_i\paren{ a_{ij}\partial_j C_{h,l}} (\bm p_{l}^{(3)}) - a C_{h,l}(\bm p_{l}^{(3)})& = \bar{f}(\bm p_{l}^{(3)}),\\
        r C_{h,l} (\bm p_{l,s}^{(1)})&= r \bar{\Phi}(\bm p_{l,s}^{(1)}),\qquad s=1,2,3,\\
        r^2 \sum_{i,j=1}^2 \beta_{ij}\partial_i C_{h,l}(\bm p_{l,s}^{(2)}) &= r^2 \bar{\Psi}(\bm p_{l,s}^{(2)}),  \qquad s=1,2.
    \end{align}
\end{subequations}

In the above formulation, we have multiplied the both side of the Neumann boundary condition and the PDE by $r$ and $r^2$, respectively. This is equivalent to a diagonal preconditioning of the linear system, such that the resulting linear system has a much better condition number. The rescaling by $r^{-1}$ in the basis and the diagonal preconditioning combined is able to make the resulting linear system not sensitive to $r$. In practice, the linear system is always solvable with condition number $\sim \mc O(10^3)$ regardless how small $r$ or $h$ is.
Notice that forming the discrete linear system \eqref{eqn:dis-cau} requires derivatives of the coefficient $a_{ij}$. For ease of implementation, we adopt a hybrid finite difference-collocation approach where the derivatives in the PDE is replaced with finite differences. 

\section{Algorithm summary}\label{sec:Algorithm sum}

To achieve a more uniform discretization of the interface curve, we begin by selecting $M$ points along the closed curve such that the arc length between adjacent points is approximately 1.5 times the grid size $h$. A cubic spline interpolation is then performed on these points to obtain a parametrized interpolant of the original interface curve, denoted by $\Gamma_h(\theta)$. Then the curve is then uniformly partitioned into $M$ discrete points with parametric coordinates based on $\theta$.

At these discrete points, different boundary integral equations \eqref{BIE for Dirichlet}, \eqref{BIE for Neumann}, \eqref{eqn:interface-bie1} corresponding to the different PDE problems can be solved by the Generalized Minimal Residual (GMRES) method. During the iteration process, all volume and boundary integrals are evaluated by solving the corresponding equivalent interface problems \eqref{eqn:equiv inter prob}, using the techniques described in sections \ref{Corrected finite difference scheme}--\ref{Local Cauchy problem}. The modified linear system \eqref{correction of scheme} can be efficiently solved using the geometric multigrid method. Once the discrete equations are solved, the grid based solution is interpolated on to the interface. After solving the boundary integral equations \eqref{BIE for Dirichlet}, \eqref{BIE for Neumann}, and \eqref{eqn:interface-bie1}, the corresponding densities can be substituted into \eqref{solution for Dirichlet}, \eqref{solution for Neumann}, and \eqref{soltion for interface problem}, respectively, to obtain the solutions to the original problems.

The main algorithmic procedure is summarized in the following pseudocode:

\begin{algorithm}[H]
\caption{Surface Interface Problem Solver}
\label{alg1}
\SetAlgoLined
\DontPrintSemicolon

\textbf{Initial Preparation}\;
1. Divide parameter domain $\Omega$ into uniform Cartesian grid\;
2. Discretize interface $\Gamma$ with equal spacing, apply cubic spline interpolation and reparameterize\;
3. Determine inside/outside status of grid points relative to interface\;

\textbf{Iterative Solution}\;
1. Formulate boundary integral equations for different problems\;
2. Convert volume integrals, single-layer and double-layer potentials to equivalent interface problems, solve using Algorithm \ref{alg2} to obtain integral values and derivatives across interface\;
3. Initialize appropriate density function, solve boundary integral equations iteratively using GMRES until convergence below preset threshold\;

\textbf{Compute the solution to the original equation} \;
1. The volume integrals, single-layer potentials, or double-layer potentials at grid points are computed using the first three steps of Algorithm \ref{alg2} under the iteratively obtained density function \;

2. Combine the values of integral according to the boundary integral expression corresponding to the original equation to obtain its solution

\end{algorithm}

\begin{algorithm}[H]
\caption{Boundary Data Computation}
\label{alg2}
\SetAlgoLined
\DontPrintSemicolon
1. Discretize interface problem using 7-point scheme, construct difference system\;
2. Solve Cauchy problem in a banded region near the interface via collocation method, add the solution as a correction to RHS of difference system\;
3. Solve linear system using multigrid method to obtain function values at grid points\;
4. Interpolate function values and derivatives across interface using nearby template points and values of correction\;
\end{algorithm}

~\\

The following analyzes the computational complexity of the algorithm. Let the discretization number of the grid be $N$. Since the uniform partitioning length of the interface is proportional to the grid discretization length, the number of interface discretization points M is approximately proportional to $N$. Therefore, in the most ideal case, the computational complexity of solving the interface problem using the multigrid method in each iteration is $\mathcal{O}(N^2)$. Additionally, during the solution of the Cauchy problem and the interpolation process, for each interface discretization point, we need to solve a small 6×6 linear system using QR decomposition, resulting in an overall complexity of $\mathcal{O}(N)$. Experimental verification shows that, when the GMRES iteration conditions remain unchanged, varying the grid density has almost no impact on the number of iteration steps. Thus, we can assume the number of convergence steps to be $\mathcal{O}(1)$. In summary, the total complexity of the algorithm is $\mathcal{O}(N^2)$.

\section{Numerical experiments}\label{sec:numer}
In this section, we present numerical results for several elliptic partial differential equations on complex bounded regions of surfaces. These equations encompass boundary value problems that include both Dirichlet and Neumann conditions, as well as interface problems on bounded regions of surfaces and interface problems on singly periodic irregular domains.

In this section, following two types of closed curves are mainly considered as the preimage of the interface $\Gamma\subset \mc{S}$: \begin{itemize}
    \item Case I: A rotated ellipse \begin{equation}\label{eqn:rot ellipse}
            \left\{\begin{array}{l}
               x=r_a \cos \theta \cos \alpha-r_b \sin \theta \sin \alpha, \\
               y=r_a \cos \theta \sin \alpha+r_b \sin \theta \cos \alpha,
              \end{array} \qquad \text { for } \theta \in[0,2 \pi),\right.
           \end{equation} where $\alpha$ is the rotation angle, $r_a$ and $r_b$ are the major and minor radius of the ellipse, respectively.\\
    \item Case II: A star-shaped domain \begin{equation}
            \left\{\begin{array}{l}
                x=\left[r_a+\epsilon \cdot \cos (m \theta)\right]\cos(\theta+\alpha), \\
                y=\left[r_b+\epsilon \cdot \cos (m \theta)\right]\sin(\theta+\alpha),
            \end{array} \qquad \text { for } \theta \in[0,2 \pi),\right.
          \end{equation} where $m$ represents the fold number and $\epsilon$ is a constant.
\end{itemize}


In each test, the GMRES iterations begin with a simple zero initial guess, and the convergence tolerance for the fixed iterations is set to $\texttt{tol} = 1.0 \times 10^{-8}$. The numerical results are presented in the form of tables and figures. Each table consists of five rows, displaying the following information: the grid size used in the Cartesian grid, the total number of discretization points on the boundary curve, the number of GMRES iterations required to solve the boundary integral equations, the maximum discrete error of the numerical solution at the interior grid nodes, and the CPU time (in seconds) on a laptop equipped with an AMD Ryzen 7 6800H processor.

\subsection{Numerical results for the boundary value problems}
\begin{example}\label{ex:bvp}
    In this example, the boundary value problems \eqref{boundary value problem} subject to both the Dirichlet boundary conditions \eqref{Dirichlet BVP} and the Neumann boundary conditions \eqref{Neumann BVP} are addressed. The underlying surface $\mc S$ is parametrized by \begin{equation}
\begin{aligned}
\bm{X}\colon & {[-1,1]^2 \rightarrow \mathbb{R}^3, } \\
& (u, v) \mapsto (3u + v, u-2v, u^3 + v^3) .
\end{aligned}\end{equation} The preimage of the interface $\Gamma$ is taken to be a rotated ellipse in Case I, with the parameters given as follows: \begin{equation}
r_a = 0.7,\qquad r_b = 0.4,\qquad \alpha = 3\pi / 5.
\end{equation} And the Dirichlet and Neumann boundary conditions on $\Gamma$ are prescribed such that the exact solution reads
\begin{equation}
   u(x, y,z) =  \mathrm{e}^{\frac{2x+y}{7}}\cos{\frac{x - 3y}{7}},\qquad (x,y,z)\in \mc{S},
\end{equation}
with $\kappa = 5$.
\end{example}

\begin{figure}[htbp]
\centering
\includegraphics[scale=0.45]{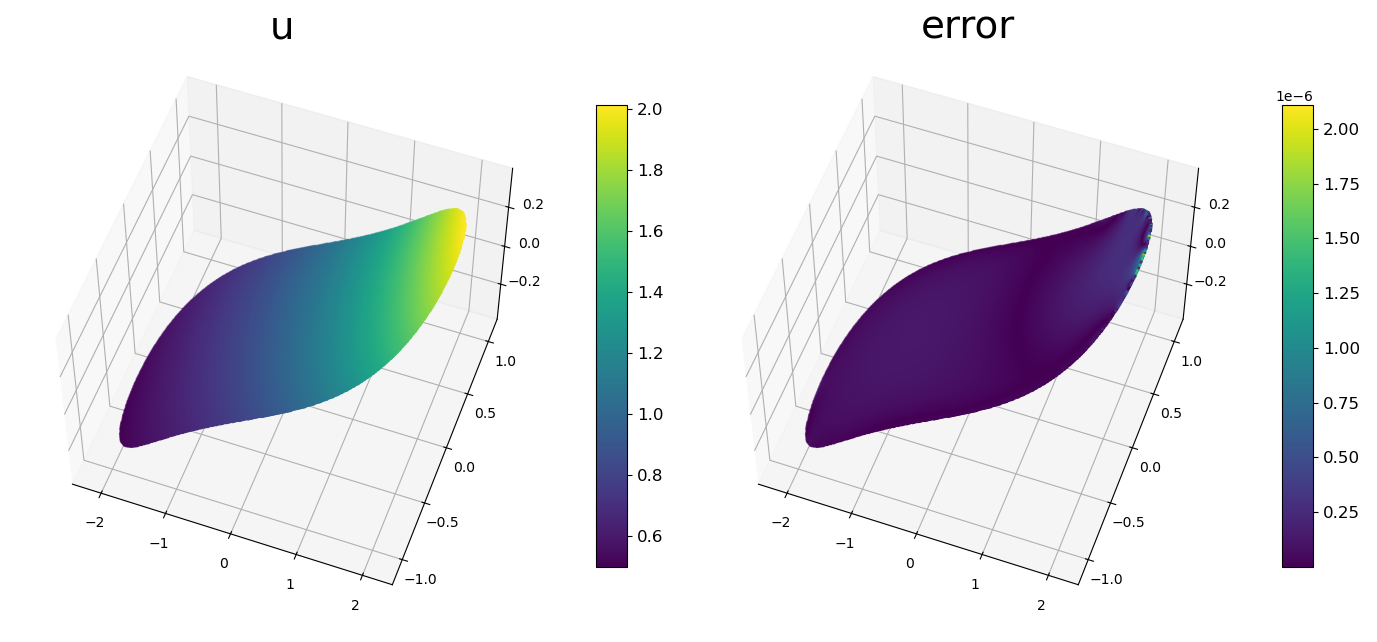}
\caption{The numerical solution and the maximum norm error in Example~\ref{ex:bvp} for the Dirichlet boundary condition on the 512 $\times$ 512 grid.}
\label{fig:Dirichlet bvp}
\end{figure}

\begin{figure}[htbp]
\centering
\includegraphics[scale=0.45]{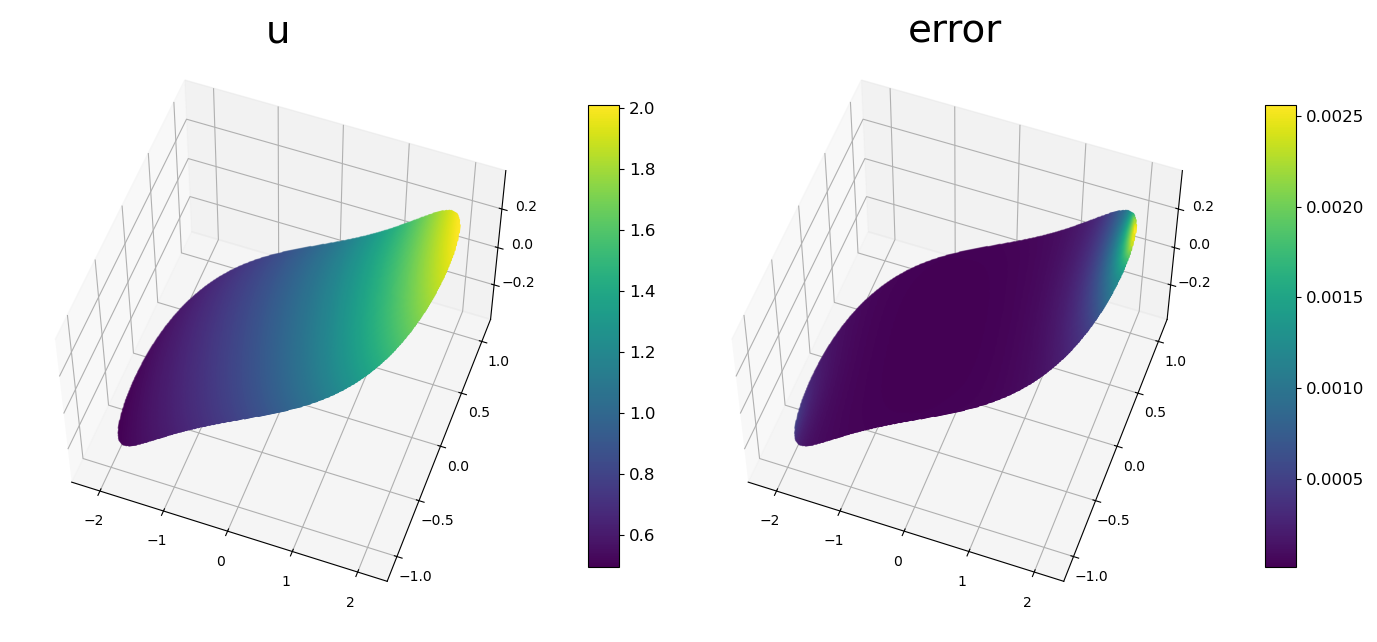}
\caption{The numerical solution and the maximum norm error in Example~\ref{ex:bvp} for the Neumann boundary condition on the 512 $\times$ 512 grid.}
\label{fig:Neumann bvp}
\end{figure}

Numerical results are summarized in Table~\ref{table:Dirichlet bvp} for the Dirichlet condition and Table~\ref{table:Neumann bvp} for the Neumann condition. Corresponding visualizations are illustrated in Fig.~\ref{fig:Dirichlet bvp} and Fig.~\ref{fig:Neumann bvp}, respectively. It can be observed that, for boundary value problems, the proposed KFBI method exhibits second-order accuracy.

\begin{table}[htbp]
\centering
\captionsetup{skip=12pt}
\renewcommand{\arraystretch}{1.25}
\begin{tabular}{|l|l|l|l|l|l|}
\hline Grid size & $64 \times 64$ & $128 \times 128$ & $256 \times 256$ & $512 \times 512$ & $1024 \times 1024$ \\
\hline M & 75 & 150 & 300 & 601 &1202  \\
\hline \#GMRES & 7 & 6 & 5 & 5 & 5 \\
\hline CPU (s) & 8.40E-1 & 3.01E+0 & 1.07E+1 & 4.47E+1 & 1.85E+2 \\
\hline $\left\|\mathbf{e}_h\right\|_{\infty}$ & 5.56E-4 & 8.97E-5 & 2.08E-5 & 2.11E-6 & 2.84E-7 \\
\hline Order & - & 2.63 & 2.11 & 3.30 & 2.89 \\
\hline
\end{tabular}
\caption{Numerical results for the Dirichlet boundary condition in Example~\ref{ex:bvp}}
\label{table:Dirichlet bvp}
\end{table}

\begin{table}[htbp]
\centering
\captionsetup{skip=12pt}
\renewcommand{\arraystretch}{1.25}
\begin{tabular}{|l|l|l|l|l|l|}
\hline Grid size & $64 \times 64$ & $128 \times 128$ & $256 \times 256$ & $512 \times 512$ & $1024 \times 1024$ \\
\hline M & 75 & 150 & 300 & 601 &1202  \\
\hline \#GMRES & 7 & 6 & 6 & 6 & 6 \\
\hline CPU (s) & 9.16E-1 & 2.89E+0 & 1.21E+1 & 4.85E+1 & 1.92E+2 \\
\hline $\left\|\mathbf{e}_h\right\|_{\infty}$ & 1.44E-1 & 3.37E-2 & 1.02E-2 & 2.43E-3 & 6.23E-4 \\
\hline Order & - & 2.10 & 1.72 & 2.07 & 1.96 \\
\hline
\end{tabular}
\caption{Numerical results for the Neumann boundary condition in Example~\ref{ex:bvp}.}
\label{table:Neumann bvp}
\end{table}

\subsection{Numerical results for the interface problems}

\begin{example}[A helicoid surface]\label{ex:inter problem helicoid}
    This example consider the interface problem \eqref{Interface problem} on a helicoid surface $\mc S$ given by \begin{equation}
\begin{aligned}
\bm{X}\colon & {[-1,1]^2 \rightarrow \mathbb{R}^3, } \\
& (u, v) \mapsto\left(u\sin{v}, u\cos{v}, v\right),
\end{aligned}
\end{equation}
where $\frac{\kappa^+}{\beta^+} = \frac{\kappa^-}{\beta^-} =1$.
A circle with radius of 0.5 is used as the preimage of the interface $\Gamma\subset\mc{S}$. The interior and exterior true solutions of the equation read
\begin{equation}
    \begin{array}{l}
        u_{i}(x, y, z) = \sin{x} \sin{y} \sin{z},   \\
        u_e(x, y, z)=(x^2 -z-1)(y^2+z-1), 
    \end{array}\qquad (x,y,z)\in \mc{S}.
\end{equation}
\end{example} 

\begin{figure}[htbp]
\centering
\includegraphics[scale=0.45]{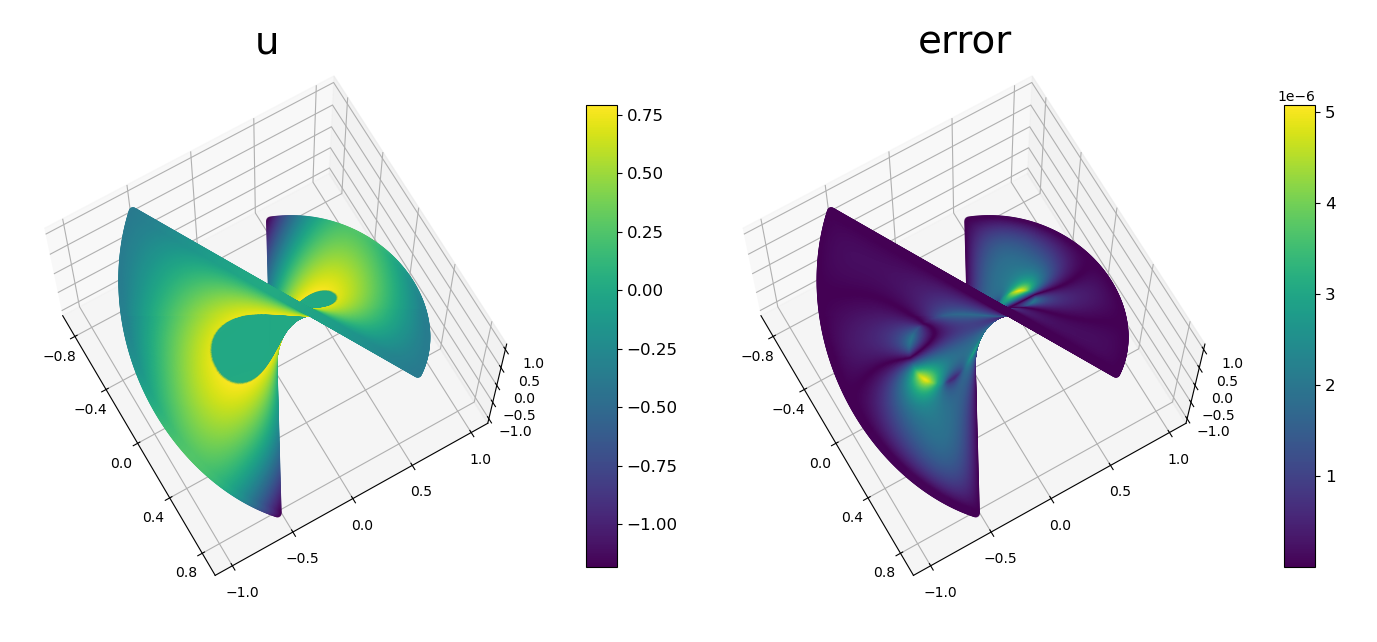}
\caption{The numerical solution and the maximum norm error in Example~\ref{ex:inter problem helicoid} on the 512 $\times$ 512 grid.}
\label{fig:inter problem helicoid}
\end{figure}

Corresponding numerical results for Example~\ref{ex:inter problem helicoid} are displayed in Table~\ref{table:inter problem helicoid} and Fig.~\ref{fig:inter problem helicoid}. The results indicate that, the proposed KFBI method achieves second-order convergence for elliptic interface problems. The low iteration counts and CPU times demonstrate the high computational efficiency of the proposed method.

\begin{table}[htbp]
\captionsetup{skip=12pt}
\renewcommand{\arraystretch}{1.25}
\centering
\begin{tabular}{|l|l|l|l|l|l|}
\hline Grid size & $64 \times 64$ & $128 \times 128$ & $256 \times 256$ & $512 \times 512$ & $1024 \times 1024$ \\
\hline M & 67 & 134 & 268 & 536 & 1072 \\
\hline \#GMRES & 7 & 6 & 6 & 6 & 6 \\
\hline CPU (s) & 9.55E-1 & 3.06E+0 & 1.24E+1 & 5.18E+1 & 2.05E+2 \\
\hline $\left\|\mathbf{e}_h\right\|_{\infty}$ & 3.48E-4 & 8.70E-5 & 2.36E-5 & 5.08E-6 & 1.23E-6 \\
\hline Order & - & 2.00 & 1.88 & 2.22 & 2.04 \\
\hline
\end{tabular}
\caption{Numerical results for the interface problem in Example~\ref{ex:inter problem helicoid}.}
\label{table:inter problem helicoid}
\end{table}

\begin{example}[A saddle surface]\label{ex:inter problem saddle}
    This example solves the interface problem \eqref{Interface problem} on a saddle surface defined by \begin{equation}
\begin{aligned}
\bm{X}\colon & {[-1,1]^2 \rightarrow \mathbb{R}^3, } \\
& (u, v) \mapsto\left(u, v, u^2-v^2\right) .
\end{aligned}
\end{equation} 
The interface $\Gamma$ is determined by a star-shaped domain in Case II with the parameters chosen as:
\begin{equation}
r_a = 0.7,\qquad r_b = 0.4,\qquad \alpha = 6\pi / 7,\qquad \epsilon = 0.3,\qquad m = 3.
\end{equation} 

The exact solutions to the equation in the interior and exterior regions are \begin{equation}
\begin{aligned}
        &u_i(x,y,z)=z\mathrm{e}^x\cos y,\qquad \\
&u_e(x,y,z)=z\mathrm{e}^y\sin x,
\end{aligned}\qquad (x,y,z)\in \mc{S}.
\end{equation} 
\end{example}

This time, we fix the grid size at $128 \times 128$. The interior-to-exterior ratios $\kappa/\beta$ are taken as a constant $c$, i.e. $\frac{\kappa^+}{\beta^+} = \frac{\kappa^-}{\beta^-} = c$. By substantially varying the value of $c$, we demonstrate that the proposed method remains applicable even under extreme ratio conditions, both in the minimal and maximal regimes. Numerical results are presented in Table~\ref{table:inter problem saddle}. 

\begin{table}[htbp]
\centering
\captionsetup{skip=12pt}
\renewcommand{\arraystretch}{1.25}
\begin{tabular}{|l|l|l|l|l|}
\hline $\kappa^+$ & 1.2E+3 & 1.2E-3 &2.0E+0 &2.0E+0 \\
\hline $\kappa^-$ & 8.0E+2 & 8.0E-4 &7.0E-1 &7.0E-1 \\
\hline $\beta^+$ & 1.2E+0 & 1.2E+0 &2.0E-3 &2.0E+3 \\
\hline $\beta^-$  & 8.0E-1 & 8.0E-1 &7.0E-4 &7.0E+2\\
\hline $c$ & 1.0E+3 & 1.0E-3 &1.0E+3 &1.0E-3\\
\hline \#GMRES& 7 & 8 &9 &13 \\
\hline $\left\|\mathbf{e}_h\right\|_{\infty}$ & 7.15E-5 & 7.40E-4  &7.22E-5 &9.97E-4\\
\hline
\end{tabular}
\caption{Numerical results for the interface problem in Example~\ref{ex:inter problem saddle} for $c = 1.0\times 10^{\pm3}$.}
\label{table:inter problem saddle}
\end{table}

Next, we fix $c = 1$. Denoting the ratio of interior-to-exterior $\beta$ as $q$, i.e. $\beta^+/\beta^- = q$. We investigate the effectiveness of the interface problem for different coefficient magnitudes by varying $q$, see Table~\ref{table:inter_problem_merged} for details.

\begin{table}[htbp]
\centering
\captionsetup{skip=10pt}
\renewcommand{\arraystretch}{1.25}
\resizebox{\textwidth}{!}{%
\begin{tabular}{|c|c|c|c|c|c|c|c|c|c|}
\hline
\multicolumn{1}{|c|}{} & \multicolumn{3}{c|}{$\beta^+ = 1$, $\beta^- = 1000$, $q = 0.001$} & \multicolumn{3}{c|}{$\beta^+ = 1$, $\beta^- = 2$, $q = 0.5$} & \multicolumn{3}{c|}{$\beta^+ = 1$, $\beta^- = 0.001$, $q = 1000$} \\
\cline{2-10}
\multicolumn{1}{|c|}{Grid size} & \#GMRES & $\|\mathbf{e}_h\|_\infty$ & Order & \#GMRES & $\|\mathbf{e}_h\|_\infty$ & Order & \#GMRES & $\|\mathbf{e}_h\|_\infty$ & Order \\
\hline
$64 \times 64$ & 16 & 1.78E-3 & - & 10 & 1.89E-3 & - & 18 & 9.44E-3 & - \\
\hline
$128 \times 128$& 16 & 5.21E-4 & 1.77 & 10 & 5.18E-4 & 1.87 & 20 & 2.63E-3 & 1.85 \\
\hline
$256 \times 256$ &16 & 1.16E-4 & 2.16 & 10 & 1.18E-4 & 2.13 & 20 & 6.33E-4 & 2.05 \\
\hline
$512 \times 512$ & 16 & 2.97E-5 & 1.97 & 10 & 3.03E-5 & 1.97 & 20 & 1.71E-4 & 1.87 \\
\hline
$1024 \times 1024$ & 16 & 7.44E-6 & 2.00 & 10 & 7.59E-6 & 2.00 & 19 & 4.49E-6 & 1.93 \\
\hline
\end{tabular}
}
\caption{Numerical results for the interface problem in Example~\ref{ex:inter problem saddle} with different coefficient ratios $q = \beta^+/\beta^-$.}
\label{table:inter_problem_merged}
\end{table}

\begin{example}[An elliptic paraboloid surface]\label{ex:inter problem elliptic Paraboloid}
    This example consider the interface problem \eqref{Interface problem} on an elliptic paraboloid surface $\mc S$ given by \begin{equation}
\begin{aligned}
\bm{X}\colon & {[-1.4,1.4]^2 \rightarrow \mathbb{R}^3, } \\
& (u, v) \mapsto\left(u, v, u^2 + v^2\right),
\end{aligned}
\end{equation}
where $\frac{\kappa^+}{\beta^+} = 3, \frac{\kappa^-}{\beta^-} =0.5$.
The interface $\Gamma$ is determined by a star-shaped domain in Case II with the parameters chosen as:
\begin{equation}
r_a = 0.7,\qquad r_b = 0.7,\qquad \alpha = 11\pi / 13,\qquad \epsilon = 0.3,\qquad m = 5.
\end{equation} The interior and exterior true solutions of the equation read
\begin{equation}
    \begin{array}{l}
        u_{i}(x, y, z) = \cos(x + y)\sin z,   \\
        u_e(x, y, z)=(x^2 -1)(y^2-1), 
    \end{array}\qquad (x,y,z)\in \mc{S}.
\end{equation}
\end{example}

\begin{figure}[htbp]
\centering
\includegraphics[scale=0.45]{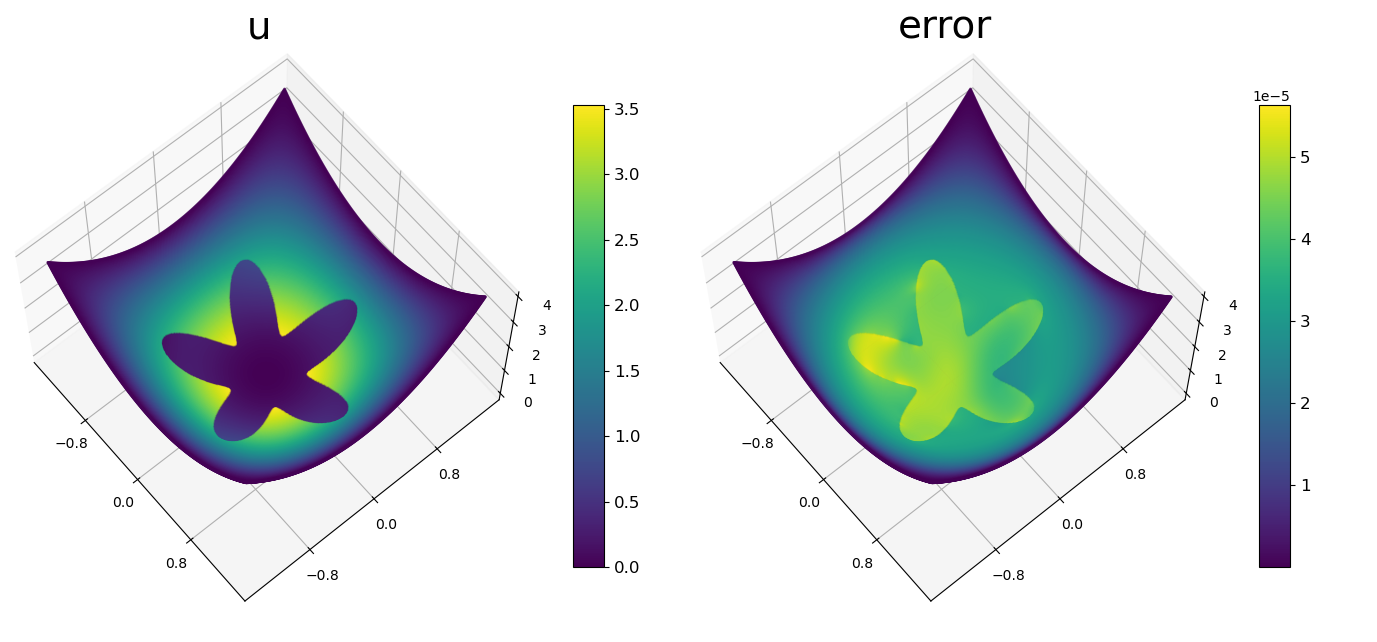}
\caption{The numerical solution and the maximum norm error in Example~\ref{ex:inter problem elliptic Paraboloid} on the 512 $\times$ 512 grid}
\label{fig:inter problem elliptic Paraboloid}
\end{figure}

Corresponding numerical results are displayed in Table~\ref{table:inter problem elliptic Paraboloid} and Fig.~\ref{fig:inter problem elliptic Paraboloid}. Example~\ref{ex:inter problem elliptic Paraboloid} considers the case where the ratio $\kappa/\beta$ differs across the interface. Numerical results show that, similar to the case with a uniform ratio, the KFBI method still exhibits second-order convergence and maintains high computational efficiency.

\begin{table}[htbp]
\captionsetup{skip=12pt}
\renewcommand{\arraystretch}{1.25}
\centering
\begin{tabular}{|l|l|l|l|l|l|}
\hline Grid size & $64 \times 64$ & $128 \times 128$ & $256 \times 256$ & $512 \times 512$ & $1024 \times 1024$ \\
\hline M & 118 & 237 &473 &946 &1892 \\
\hline \#GMRES & 10 & 10 & 9 & 10 & 10 \\
\hline CPU (s) & 7.09E+0 & 2.77E+1 & 1.00E+2 & 4.37E+2 & 1.99E+3 \\
\hline $\left\|\mathbf{e}_h\right\|_{\infty}$ & 4.10E-3 & 8.66E-4 & 1.96E-4 & 5.63E-5 & 1.38E-5 \\
\hline Order & - & 2.24 & 2.15 & 1.80 & 2.03\\
\hline
\end{tabular}
\caption{Numerical results for the interface problem in Example~\ref{ex:inter problem elliptic Paraboloid}.}
\label{table:inter problem elliptic Paraboloid}
\end{table}

\begin{figure}[htbp]
\centering
\includegraphics[scale=0.6]{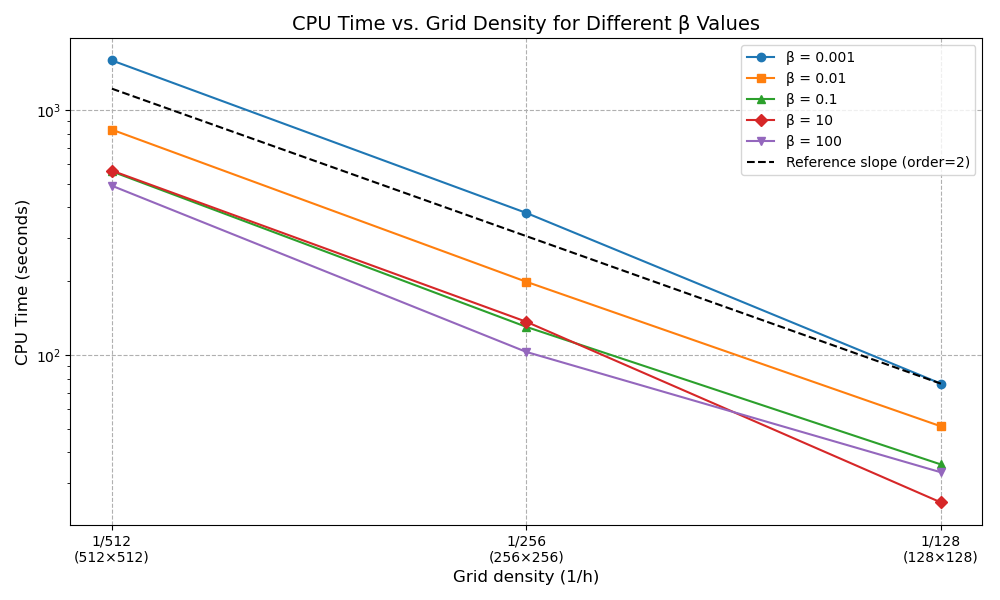}
\caption{The CPU time in Example~\ref{ex:inter problem elliptic Paraboloid} with different interior-exterior ratios.}
\label{fig:the CPU time inter problem elliptic Paraboloid with different interior-exterior ratios}
\end{figure}

We now fix both interior/exterior $\kappa$ and interior $\beta$ at 1, while varying the exterior $\beta$ to investigate how different interior-exterior ratios affect the computational results. As demonstrated in previous test cases, the number of GMRES iterations remains virtually unaffected by grid refinement. Therefore, iteration counts will not be discussed here. The corresponding maximum norm errors and CPU time consumption are presented in the following Fig.~\ref{fig:the error of inter problem elliptic Paraboloid with different interior-exterior ratios}--Fig.~\ref{fig:the CPU time inter problem elliptic Paraboloid with different interior-exterior ratios}.

\begin{figure}[htbp]
\centering
\includegraphics[scale=0.6]{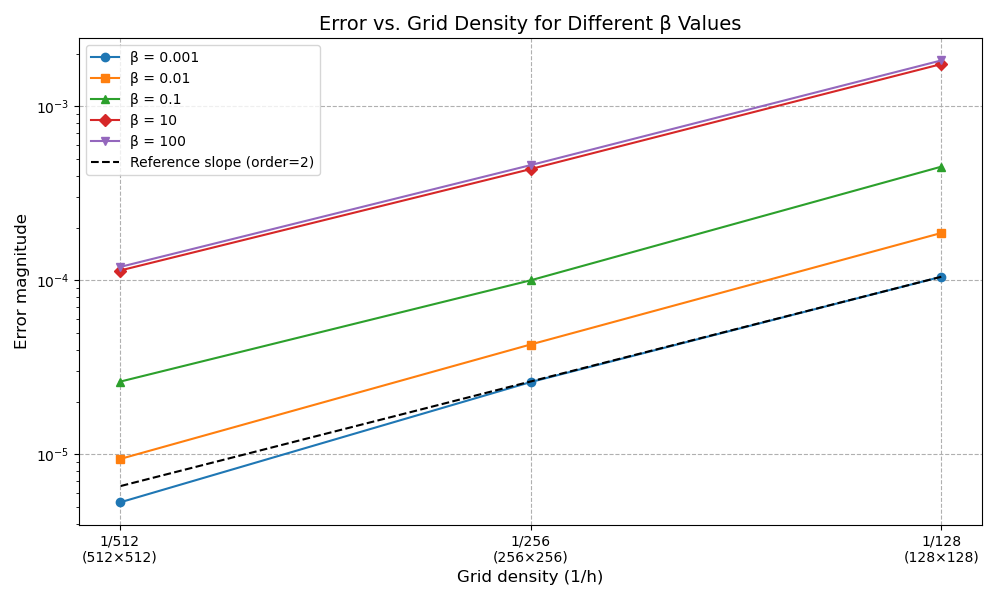}
\caption{The maximum norm error in Example~\ref{ex:inter problem elliptic Paraboloid} with different interior-exterior ratios.}
\label{fig:the error of inter problem elliptic Paraboloid with different interior-exterior ratios}
\end{figure}

\subsection{Interface problems on closed surfaces}

\begin{example}[A torus]  \label{ex:inter problem torus}
In this example, we consider the interface problem on a torus $\mathbb{T}^2\subset\mathbb{R}^3$ with a global parametrization given by \begin{equation}
        \begin{aligned}
            \bm{X}\colon &(\mathbb{R}/2\pi\mathbb{Z})^2\to\mathbb{R}^3\\
            &(u,v)\mapsto\left((R+r\sin u)\cos v, (R+r\sin u)\sin v,r\cos u\right),
        \end{aligned}
    \end{equation} where $R,r>0$ are the major and minor radii, respectively. 

    Parameters are taken as follows: \begin{subequations}
        \begin{align}
            &R=2.0,\qquad\qquad\qquad r=0.8,\\
            &\beta^\pm=1.25\pm0.75,\qquad\kappa^\pm=1.25\pm 0.75.
        \end{align}
    \end{subequations} Two types of surfaces $\Gamma_1$ and $\Gamma_2$ in Case I and Case II are selected as test-case interfaces, with their respective parameter chosen as:
\begin{equation}
\begin{aligned}
&\Gamma_1\colon\quad r_a = 1.0,\qquad r_b = 0.6,\qquad \alpha = 9\pi / 13,\\
&\Gamma_2\colon \quad  r_a = 0.6,\qquad r_b = 0.6,\qquad \alpha = \pi / 4,\qquad \epsilon = 0.4,\qquad m = 3.
\end{aligned}
\end{equation} 
The interior and exterior true solutions of the equation read
\begin{equation}
    \begin{array}{l}
        u_i(u, v) = \sin{u} \cos{v},   \\
        u_e(u , v)= \cos{u} \sin{v}, 
    \end{array}\qquad (u, v)\in \Omega.
\end{equation}
\end{example}

\begin{figure}[htbp]
\centering
\includegraphics[scale=0.45]{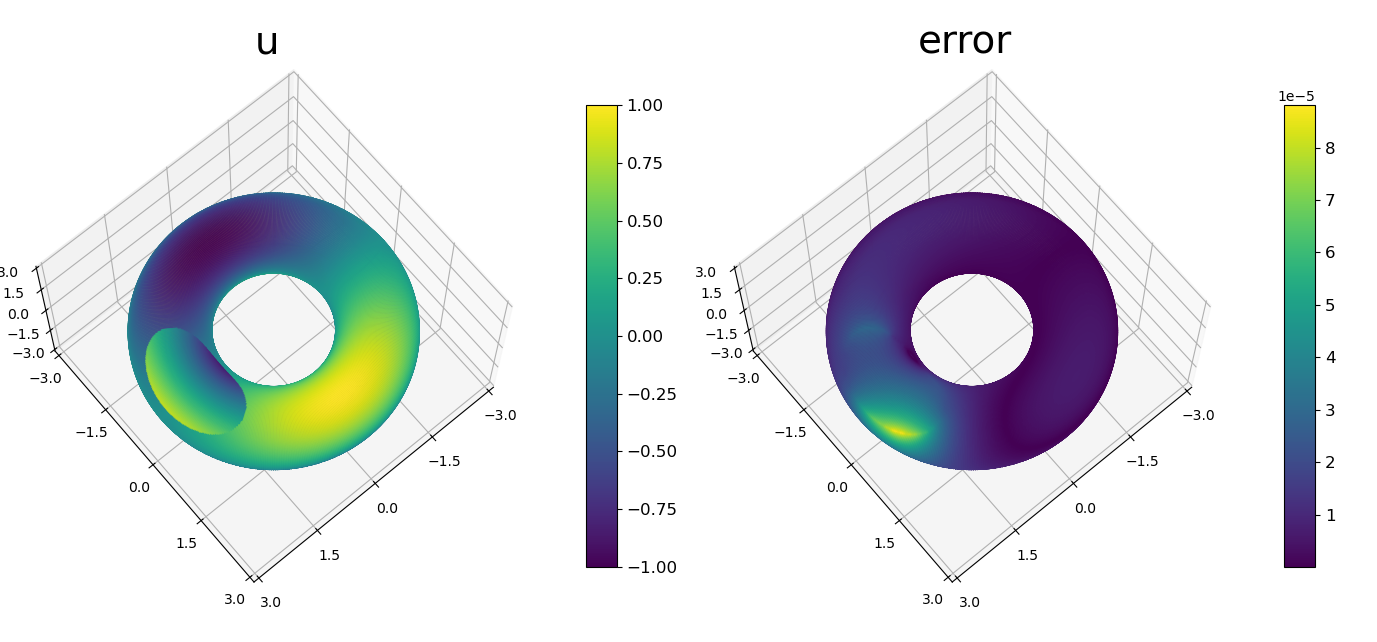}
\caption{The numerical solution and the maximum norm error in Example~\ref{ex:inter problem elliptic Paraboloid} of $\Gamma_1$ on the 512 $\times$ 512 grid.}
\label{fig:inter problem torus2}
\end{figure}

\begin{figure}[htbp]
\centering
\includegraphics[scale=0.45]{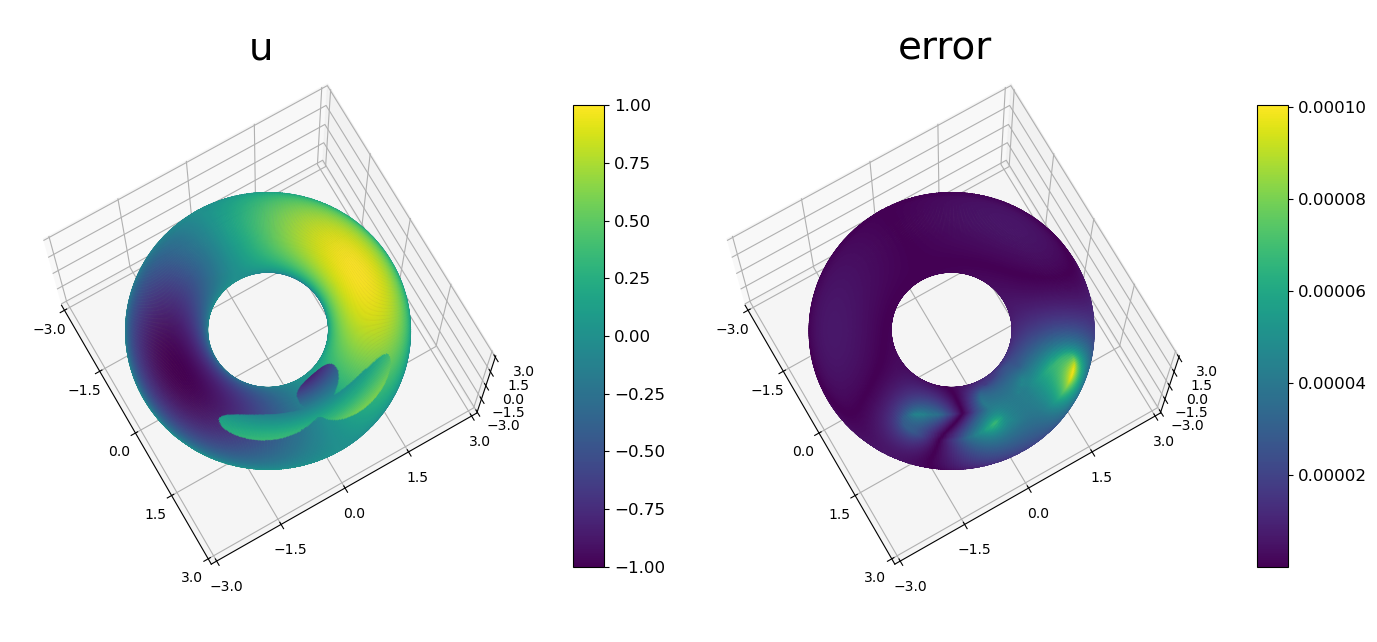}
\caption{The numerical solution and the maximum norm error in Example~\ref{ex:inter problem elliptic Paraboloid} of $\Gamma_2$ on the 512 $\times$ 512 grid.}
\label{fig:inter problem torus}
\end{figure}

The computational results and maximum norm errors are shown in Fig.~\ref{fig:inter problem torus2} and Fig.~\ref{fig:inter problem torus}. The numerical performance for $\Gamma_1$ under progressive grid refinement is summarized in Table \ref{table:inter problem torus}, including iteration numbers, computational time, and final errors. The results demonstrate that the method achieves second-order convergence accuracy.

\begin{table}[htbp]
\centering
\captionsetup{skip=12pt}
\renewcommand{\arraystretch}{1.25}
\begin{tabular}{|l|l|l|l|l|l|}
\hline Grid size & $64 \times 64$ & $128 \times 128$ & $256 \times 256$ & $512 \times 512$ & $1024 \times 1024$ \\
\hline M & 44 & 87 & 175 &350 &700  \\
\hline \#GMRES & 17 &17 &17 &14 &15 \\
\hline CPU (s) & 3.70E+0 & 1.53E+1 & 6.20E+1 & 2.13E+2 & 1.06E+3 \\
\hline $\left\|\mathbf{e}_h\right\|_{\infty}$ & 1.37E-2 & 2.12E-3 & 3.56E-4 & 1.00E-4 & 2.09E-5 \\
\hline Order & - & 2.69 & 2.58 &1.83 & 2.26 \\
\hline
\end{tabular}
\caption{Numerical results for the interface problem in Example~\ref{ex:inter problem torus} of $\Gamma_1$.}
\label{table:inter problem torus}
\end{table}

\begin{example}[A Dupin cyclide] \label{ex:inter problem Dupin cyclide}
This example considers the interface problem on a Dupin cyclide $\mathcal{D}\subset\mathbb{R}^3$, whose parametrization is given by
    \begin{equation}
        \begin{aligned}
            \bm{X}\colon &(\mathbb{R}/2\pi\mathbb{Z})^2\to\mathbb{R}^3\\
            &(u,v)\mapsto\left(\frac{d(c-a\cos u\cos v)+b^2\cos u}{a-c\cos u\cos v},\frac{b\sin u\,(a-d\cos v)}{a-c\cos u\cos v},\frac{b\sin v\,(c\cos u-d)}{a-c\cos u\cos v}\right).
        \end{aligned}
    \end{equation} The parameters are chosen as $a=1,b=1,c=-0.3,d=0.5$.

The diffusion and reaction coefficients are set to be $\beta^\pm = 1.9\pm 1.1$ and $\kappa^\pm=1.9\pm 1.1$. The preimage of the interface in parameter space is a unit circle, and we set the interior and exterior exact solutions to the equation as
\begin{equation}
    \begin{array}{l}
        u_i(u, v) = \sin{u} \sin{v},   \\
        u_e(u , v)= \cos{u} \cos{v}, 
    \end{array}\qquad (u, v)\in \Omega.
\end{equation}
\end{example}

\begin{figure}[htbp]
\centering
\includegraphics[scale=0.45]{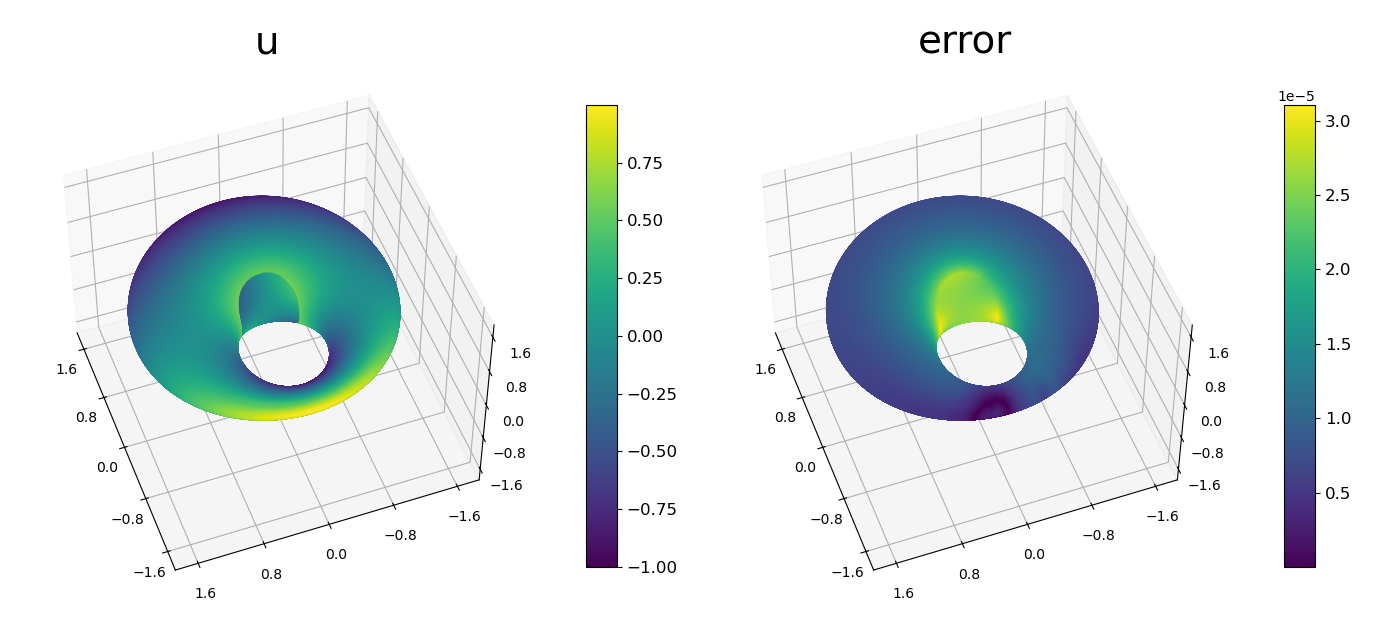}
\caption{The numerical solution and the maximum norm error in Example~\ref{ex:inter problem Dupin cyclide} of $\Gamma_2$ on the 1024 $\times$ 1024 grid.}
\label{fig:inter problem Dupin cyclide}
\end{figure}

The numerical results and visualizations can be found in Table~\ref{table:inter problem Dupin} and Fig.~\ref{fig:inter problem Dupin cyclide}.

\begin{table}[htbp]
\centering
\captionsetup{skip=12pt}
\renewcommand{\arraystretch}{1.25}
\begin{tabular}{|l|l|l|l|l|l|}
\hline Grid size & $64 \times 64$ & $128 \times 128$ & $256 \times 256$ & $512 \times 512$ & $1024 \times 1024$ \\
\hline M & 43 & 85 & 171 &341 &683  \\
\hline \#GMRES & 8 &8 &8 &7 &7 \\
\hline CPU (s) & 6.43E+0 & 2.55E+1 & 9.97E+1 & 3.22E+2 & 1.39E+3 \\
\hline $\left\|\mathbf{e}_h\right\|_{\infty}$ & 7.44E-3 & 1.95E-3 & 5.03E-4 & 1.22E-4 & 3.10E-5 \\
\hline Order & - & 1.93 & 1.96 &2.04 &1.97 \\
\hline
\end{tabular}
\caption{Numerical results for the interface problem in Example~\ref{ex:inter problem Dupin cyclide}.}
\label{table:inter problem Dupin}
\end{table}

\section{Conclusion}\label{sec:conclusion}

This study develops an efficient kernel-free boundary integral method for accurately solving elliptic boundary value and interface problems on surfaces. The proposed method reformulates the original problem as boundary integral equations, and evaluates the resulting singular integrals indirectly by solving equivalent interface problems on Cartesian grids via a geometric multigrid method. This strategy avoids the difficulties associated with mesh transformation and direct numerical quadrature. The combination of boundary integral equations with finite difference methods enables the proposed approach to effectively handle complex geometric interfaces and discontinuities near the interface, offering broad applicability.

The proposed method offers several key advantages, including flexible treatment of complex geometries, high computational efficiency, robust convergence rates, and mesh-independent convergence behavior of the iterative solvers. Numerical results demonstrate that the method maintains second-order accuracy regardless of whether the interior/exterior $\kappa$ and $\beta$ ratios are equal or divergent, while showing robust applicability across wide $\beta$ variations. For interface problems on closed surfaces with periodic boundary conditions, the proposed method maintains second-order convergence and exhibits favorable computational efficiency.

The current method is confined to surfaces that possess a regular global parameterization. Nevertheless, such parameterization is not anticipated for general closed surfaces. Despite this limitation, the current potential theory-based approach offers a framework for efficiently solving interface problems on surfaces through an immersed interface-type technique. To extend the applicability of this method to a wider range of closed surfaces, we will employ alternative surface and interface discretization techniques and numerical schemes, such as those described in \cite{ying2013kernel} and \cite{beale2020solving}.


\section*{CRediT authorship contribution statement}

Pensong Yin: Conceptualization, Methodology, Software, Investigation, Writing - Original Draft.\\

Wenjun Ying: Supervision, Writing - Review \& Editing.\\

Yulin Zhang: Project administration, Supervision, Writing - Review \& Editing.\\

Han Zhou:  Methodology, Project administration, Supervision, Writing - Review \& Editing.\\


\section*{Declaration of competing interest}

The authors declare that they have no known competing financial interests or personal relationships that could have appeared to influence the work reported in this paper.

\section*{Data availability}

Data will be made available on request.

\section*{Acknowledgement}

W. Ying is supported by the Shanghai Science and Technology Innovation Action Plan in Basic Research Area (Project No. 22JC1401700), the National Natural Science Foundation of China in the Division of Mathematical Sciences (Project No. 12471342). It is also partially supported by the National Key R\&D Program of China (Project No. 2020YFA0712000) and the fundamental research funds for the central universities.

\end{document}